\input amstex
%&AMS-TeX
\magnification=1200
\input amsppt.sty

\parindent 20 pt
 
\NoBlackBoxes
\define\df{\dsize\frac}
\define\bk{\bigskip}
 \define\mk{\medskip} 

\define \be{\beta}

\define \g{\gamma}

\define\ri{\rightarrow}
\define\Ri{\Rightarrow}
\define \s{\sigma}
\define\si{\sigma}
\define \fa{\forall}

\define \tY{\tilde{Y}}\define \tg{\tilde{g}}\define \tT{\tilde{T}}
\define \tZ{\tilde{Z}}
\define \tB{\tilde{B}}\define \tP{\tilde{P}}
\define \tX{\tilde{X}}\define \tb{\tilde{b}}
\define \la{\langle}
\define \ra{\rangle}

\define \G{\Gamma}

\define \dl{\delta}

\define \1{^{-1}}
\define \2{^{-2}}
\define \p{\partial}

\define \Aut{\operatorname{Aut}}

\define \Ab{\operatorname{Ab}}
\define \Center{\operatorname{Center}}
\define \un{\underline}

\define \ub{\underbar}
\define \pf{\demo{Proof}}
\define \fc{\frac}

\define \edm{\enddemo}
\define \ep{\endproclaim}

\baselineskip 20pt
\topmatter
\heading{\bf On the quotient of the braid group by commutators of
transversal
half-twists and its group actions}\rm\endheading 
\mk

\centerline{M. Teicher}
\abstract\nofrills{\bf Abstract.}\ This paper presents and describes a
quotient of the Artin
braid group by commutators of transversal half-twists and investigates
its group actions.
We denote the quotient by $\tilde B_n$ and refer to the groups which
admit an action of $\tilde B_n$ as $\tilde B_n$-groups.
The group $\tilde B_n$ is an extension  of a solvable group by a
symmetric group. We 
 distinguish special elements in
$\tilde B_n$-groups which we call prime elements and we  give a
criterion for an element
to be prime.
$\tilde B_n$-groups appear as fundamental groups of complements of
branch curves.
\endabstract
\subjclass
20F36\endsubjclass\endtopmatter
 
In this paper we describe a quotient of the Artin braid group by
commutators of
transversal half-twists and investigate its group actions.
These groups  turned out to be extremely important in describing
fundamental
groups of complements of branch curves (see, e.g.,   \cite{Te}).
The description here is completely independent from the
algebraic-geometrical background and
provides an algebraic study of the groups involved, using a topological
approach to the braid
group.
 We denote the quotient by $\tilde B_n$ and refer to the groups which
admit an action of
$\tilde B_n$ as $\tilde B_n$-groups.
 In particular, we study $\tilde P_n$, the image of the pure braid group
in $\tilde B_n,$ and
prove that $\tilde B_n$ is an extension of a solvable group by a
symmetric group.  
The main results on the structure of $\tilde B_n$ are Theorem 6.4 and
Corollary 6.5.
 We  
distinguish special elements in $\tilde B_n$-groups which we  call 
prime elements, compute
the action of half-twists on prime elements (\S2-\S4), and finally we 
give a criterion for an
element to be prime (see Proposition 7.1). This criterion will be
applied to the study of
fundamental groups of complements of branch curves.

Fundamental groups related to algebraic varieties are very important in
classification
problems and in topological studies in algebraic geometry.
These groups are very difficult to compute.
The group $\tilde B_n$ and the $\tilde B_n$-groups appeared when we were
computing such groups.
It turns out that {\it all} new examples of fundamental groups of
complements of branch curves
are $\tilde B_n$-groups and unlike previous expectations they are
``almost solvable'' (like
$\tilde B_n$ itself) . Our future plans are to study $\tilde B_n$-groups
indpendent of the
fundamental group problems. The ultimate goal is to classify all $\tilde
B_n$-groups.

The paper is divided as follows:
\roster
%\item"0." Introduction
\item"1." Definition of $\tilde B_n$
\item"2." $\tB_n$-groups and prime elements
\item"3." Polarized pairs and uniqueness of coherent pairs
\item"4." $\tB_n$-action of  half-twists
\item"5." Commutativity properties
\item"6." On the structure of $\tilde B_n$ and  $\tilde P_n $  
\item"7."  Criterion for prime element
\endroster

Throughout this paper we use the following notations:
 $D$ is a disc, $K$ is a finite subset of $ D,$ \ $n=\# K$\ $(n\ge 4),$
 \ the braid group $B_n=B_n[D,K]$  is the group  of all  diffeomorphisms
$\be: D
\to D$ which preserve $K$ and act as an identity on $\partial D,$ under
the equivalence
relation that $\be_1\sim \be_2$ if their actions on $\pi_1(D-K,*)$
coincide.
For any element $X,Y$ in a group $G$ we write  $X_Y=Y^{-1}XY.$
 (The ordinary notation for conjugation is $X^Y=YXY^{-1}.)$
\bk
\subheading{\S1. Definition of $\bold{\tilde B_n}$}

In this section we define the group $\tilde B_n.$
This group and the groups on which it acts (called  $\tilde B_n$-groups)
are the central
objects of our investigation.
We also introduce the basic notions of a frame and a good quadrangle.
In Claim 1.1 and Lemma 1.2 we establish certain basic identities which
will be used
repeatedly in later sections.

We recall here the definition of a half-twist in the braid group.

\newpage

\definition{Definition} \ $\un{\text{Half-twist w.r.t.}\
\left[-\df{1}{2},\frac{1}{2}\right] }$

Consider $D_1,$ the  unit disc, $\pm\df{1}{2}\in D_1.$
Take $\rho: [0,1]\to[0,1]$ continuous such that $\rho(r)=\pi$ for 
$r\le\df{1}{2}$ and $\rho(1)=0.$
Define $\dl: D_1\to D_1: \dl(re^{i\theta})=re^{i(\theta+\rho(r))}.$
Clearly, $\dl\left(\df{1}{2}\right)=-\df{1}{2},$\quad
$\dl\left(-\df{1}{2}\right)=\df{1}{2},$ and $\dl|_{\p D_{1}}=Id.$
The disc of radius $\df{1}{2}$ rotates $180^\circ$ counterclockwise.
Outside of this disc it rotates by smaller and smaller angles till it is
fixed
on the unit circle.
Thus we get a braid $[\dl]\in
B_2\left[D_1,\left\{\pm\df{1}{2}\right\}\right].$ $[\dl]$ is called the
half-twist w.r.t. the segment
$\left[-\df{1}{2},\df{1}{2}\right].$\enddefinition \bigskip
Using the above definition we define a general
half-twist.
\definition{Definition} \ $\un{H(\si), \ \text{half-twist w.r.t. a
path}\ \sigma}.$

 Let $D,\ K$ be as above, $a,b\in K.$
Let $\sigma$ be a path from $a$ to $b$ which does not meet any other
point of $K.$   We take  a small topological disc $D_2\subset D$ such
that
$\sigma\subset D_2,\quad D_2\cap K=\{a,b\}.$
We take a diffeomorphism $\psi: D_2\to D_1$ (unit disc) such that 
$\psi(\sigma)=\left[-\fc{1}{2},\fc{1}{2}\right].$
\ $\psi(a)=-\fc{1}{2}$\quad $\psi(b)=\fc{1}{2}.$
 We consider a ``rotation'' $\psi\dl\psi\1: D_2\to D_2;$
 $\psi\dl\psi\1$ is the identity on the boundary of $D_2.$
We extend it to $D$ by the identity.
We define the half-twist $H(\sigma)$ to be the conjugacy  class of the
extension of
$\psi\dl\psi\1.$ \enddefinition

\midspace{2.20in}\caption{Fig. 1.0}

We sometimes denote the path by $x$ and the corresponding half-twist
$H(x)$ by $X.$

\definition{Definition}\ Let $D,K$ be as denoted above.
Let $H(\s_1)$ and $H(\s_2)$ be 2 half-twists in $B_n=B_n[D,K].$
We say that $H(\s_1)$ and $H(\s_2)$ are:

\roster
\item"(i)" \ub{weakly disjoint} \ if  $\s_1\cap\s_2\cap K=\emptyset.$
\item"(ii)"  \ub{transversal}  \ if $\s_1$ and $\s_2$ are weakly
disjoint
and  intersect each other exactly
once (and not in any point of $K),$ i.e.,  
$\s_1\cap\s_2=$ one point,\ $\s_1\cap\s_2\cap K=\emptyset.$
\item"(iii)"  \ub{disjoint} \ if $\s_1\cap\s_2=\emptyset.
$
\item"(iv)" \ub{adjacent} \  if  $\s_1\cap\s_2\cap K=$ one point.
\item"(v)" \ub{consecutive} \  if they are adjacent and $\s_1\cap \s_2$
do
not intersect outside of $K,$
i.e.,  $\s_1\cap\s_2=$ point $\in K.$ 
\item"(vi)" \ub{cyclic} \quad if $\s_1\cap \s_2=$ 2 points $\in K.$
\endroster

\enddefinition

\midspace{2.70in}\caption{Fig. 1.1}

\proclaim{Claim  1.0}
Let $X,Y$ be two half-twists in $B_n.$ Then:
\roster  
\item"(i)"If $X,Y$ are disjoint, then $[X,Y]=1,$ i.e., $XY=YX.$
\item"(ii)" If $X,Y$ are consecutive, then $\la X,Y\ra=
XYXY\1X\1Y\1=1,$ thus\newline $XYX=YXY$,\
$X_{Y\1}=Y_X$ and $X_{Y\1X\1}=Y.$\newline  We say then that $X$ and $Y$
satisfy the triple
relation. \item"(iii)" If $X=H(x) $ is represented by a diffeomorphism
$\be$ and 
$Y=H(y),$ then  $Y_X = X\1YX=H((y)\be).$
If $X$ and $Y$ are consecutive, then $(y)\be$, $x$ and $y$ create  a
triangle whose interior
contains no points of $K$ (see Fig. 1.2).
\item"(iv)" If $X_1,X_2,X_3$ are 3 consecutive half-twists (see for
example Fig.
1.4) then
$[X_2,(X_2)_{X_1X_3}]=(X_3)^{-2}_{X_2\1}(X_1)^{-2}_{X_2\1}X_1^2X_3^2.$\endroster\ep

\midspace{2.00in}\caption{Fig.   1.2}

\demo{Proof} We will use the following observation: If a half-twist in
$B_n$ is determined by a $180^\circ$ rotation around
$\sigma$, then its inverse is determined by a clockwise rotation.
An equation in $B_n[D,K]$ may always be verified by studying the action
of each side of the
equation on a set of generators of the fundamental group $
\pi_1(D-K,*).$
The action of a half-twist on a simple loop is easy to determine:
If $a,b
\in K,$ and $\G_1$ and $\G_2$ are 2 simple loops  around $a$ and $b,$
respectively, and if $
\sigma$ connects $a$ and $b$ \ $(\sigma\cap K=\{a,b\})$, then
$H(\sigma)$ transfers $\G_1$ to
$\G_2$, $\G_2$ to $\G_2\G_1\G_2\1$, and does not move any simple loop
around any other point
of $K.$ $H(\sigma)\1$ transfers $\G_1$ to $\G_1\1\G_2\G_1$ and $\G_2$ to
$\G_1.$

Claims (i), (ii) and (iii) follow from  the above abservation.
For example, for (iii), and $X,Y$ consecutive half-twists  with $x\cap
y  = b$ note that
$\be$ transfers $b$ to $c$ and fixes $a$ (see Fig. 1.2).
Thus it transfers $y$ to a path connecting $a$ and $c$.
Thus $H((y)\be)$ fixes $b$ and interchanges $a$ and $c.$
The same is true for $X\1 YX$ and thus they coincide.
 If $X$ and $Y$ are disjoint, (iii) is
a   trivial result of (i) and the above observation.
Items (i) and (ii) are well-known properties
of the braid group.

Claim (iv) is a consequence of (i) and (ii) as follows:

In $B_n:[X_1,X_3]=1$ and $\la X_1,X_2\ra=\la X_3, X_2\ra=1.$
Thus by (ii), $X_i X_j\1 X_i\1=X_j\1 X_i\1 X_j$ for $(i,j)=(2,3)$ or
(1,2) or (3,2) or
(2,1).  We shall use these relations in the following list of equations:

$\qquad\qquad\quad[ X_2,( X_2) _{ X_1 X_3}
]$

$\qquad\qquad\quad \quad =[ X_2, X_3\1  X_1\1 X_2 X_1 X_3]$
\smallskip
$\qquad\qquad\quad \ \ \ =  X_2 X_3\1
 X_1\1 X_2\tX_1{\underbrace{ X_3 X_2\1 X_3\1}} X_1\1 X_2\1 X_1 X_3$
\baselineskip .03pt
\flushpar $\qquad  \qquad\qquad \qquad \qquad\qquad
\qquad\quad\qquad\quad\  
\overbrace{\ \quad \  \qquad}$ 

$\text{(Claim\ 1.0(ii))} =  X_2X_3\1
X_1\1X_2X_1X_2\1X_3\1{\underbrace{X_2X_1\1X_2\1}}X_1X_3$  

$\text{(Claim\ 1.0(ii))}   \qquad \qquad
\qquad \qquad \qquad \qquad \ \quad {\overbrace{X_1\1X_2\1X_1}}$
\baselineskip .03pt
\flushpar $\ \qquad \qquad \qquad \qquad\qquad \qquad \qquad\quad
\qquad \qquad \  \underbrace{\ \quad \  \qquad}$ 
\flushpar $\qquad\ \qquad \qquad \qquad\qquad \qquad \qquad\quad\
\qquad \qquad \  \overbrace{\ \quad \  \qquad}$ 
\flushpar $\qquad\   ([X_1,X_3]=1)  =\  X_2X_3\1X_1\1X_2
\underbrace{X_1X_2\1X_1\1}\ X_3\1X_2\1X_1X_1X_3$
\baselineskip 10pt
\flushpar $\qquad\qquad \qquad\qquad   \quad \qquad \qquad
\qquad  \quad \ \ {\overbrace{X_2\1X_1\1X_2}}$
\baselineskip 20pt

$\qquad\qquad\qquad = \quad \ 
X_2X_3\1X_1^{-2}\ {\underbrace{X_2\ X_3\1\ X_2\1}}X_1^2X_3$

$ \text{(Claim\ 1.0(ii))}  \qquad\qquad\qquad\quad\ \ 
{\overbrace{X_3\1X_2\1X_3}}$

$\qquad\qquad\qquad =\  X_2X_3^{-2}X_1^{-2}X_2\1
X_1^2X_3^2$

$\qquad\qquad\qquad =\ (X_3)^{-2}_{X_{2}\1}\ (X_1)^{-2}_{X_{2}\1}
X_1^2 X_3^2.$

\qquad \qed\enddemo 

\remark{Remark}

(a) $XYX=YXY$ is called the {\it triple relation}.

	(b) Claim (iv) will be used in the proof of Lemma 1.2 which is the
first step in
understanding $\tilde B_n.$\endremark

\newpage
 
\definition{Definition} \ $\un{\tilde B_n}$

Let $n\ge 4.$  $\tilde{B}_n$ is  the quotient of $B_n,$ the braid group
of order $n,$
by the subgroup normally generated by the commutators $[H(\s_1),
H(\s_2)]$, where
$H(\s_1)$ and $H(\s_2)$ are transversal half-twists.\enddefinition

\demo{Notation}\

 Let $Y\in B_n.$  We denote the image of $Y$ in $\tilde{B}_n$
by  $\tilde{Y}.$
If $Y$ is a half-twist in $B_n$ we call $\tilde{Y}$  a half-twist in
$\tilde{B}_n.$
We call two half-twists $\tilde{Y}, \tilde{X}$ in $\tilde{B}_n$ disjoint
(or
weakly disjoint, adjacent, consecutive, transversal) if $Y,X$ are
disjoint
(or
weakly disjoint, adjacent, consecutive, transversal).

\edm
\definition{Definition}\  $\un{\text{A frame of $\tilde B_n$}}$

Let $K=\{a_1\dots a_n\}$.
Let $\sigma_i$,\ $i=1
\dots n-1,$ be a consecutive sequence of simple paths in $D$ connecting
the points of $K$ such
that $\sigma_i\cap K=\{a_i, a_{i+1}\}$ and
$$\sigma_i\cap\sigma_j=\cases\emptyset\  & i\ne
j,j+1\\ a_j\  & i=j-1\\ a_{j+1}\  & i=j+1\endcases,$$  $
\{X_i=H(\sigma_i)\}$ is a frame of
$B_n$ and $\{\tilde X_i\}$ is a frame of $\tilde B_n.$
\enddefinition
\remark{Remark}
We also refer to a frame  as a standard base of $\tB_n.$
By the classical Artin theorem, a frame generates $\tilde B_n$ (see
\cite{A}).\endremark

\definition{Definition}\  $\un{\text{Polarized half-twist,
polarization}}$

We say that a half-twist $X\in B_n$ (or $\tX$ in $\tilde B_n$) is
polarized if
we choose an order on the end points of $X.$ The order is called the
polarization of $X$ (or $\tX$).\enddefinition

\definition{Definition }\ \underbar{Orderly adjacent}

Let $X,Y$ (or $\tilde X,\tilde Y)$  be two adjacent polarized
half-twists in $B_n$ (resp.
in $\tilde{B}_n$).  We say that $X,Y$ (or $\tilde X,\tilde Y)$  are 
{\it orderly adjacent}
if their common point is the ``end'' of one of them and the
``origin'' of another.
\enddefinition

\newpage

\definition{Definition}\ \underbar{Good quadrangle}

Let $H(\sigma_i),$\ $i=1,\dots,4,$ be four half-twists such that
$H(\sigma_i)$ is
consecutive to $H(\sigma_{i+1})\pmod 4,$\ $H(\sigma_i)$ is disjoint from
$H(\sigma_{i+2})\pmod 4$ and in the disc with boundary
$\bigcup\limits_{i=1}^4\sigma_i$
there are no points of $K$ (see Fig. 1.3).
 
 We say that $\{H(\s_i)\}$ is a good
quadrangle in $B_n$, and $\{\widetilde{H(\s}_i)\}$ is a good quadrangle
in
$\tilde{B}_n.$\enddefinition

\midspace{2.00in}\caption{Fig.   1.3}

\proclaim{Claim   1.1}
\roster\item"(a)" transversal, disjoint $\Rightarrow$ weakly disjoint.
\newline consecutive $\Rightarrow$ adjacent.
\item"(b)"  Every two transversal or disjoint half-twists   commute.
\newline Every two consecutive half-twists   satisfy the
triple relation $(\tilde X\tilde Y\tilde X=\tilde Y\tilde X\tilde Y$ or
$\tilde X\1\tilde
Y\tilde X=\tilde Y\tilde X\tilde Y\1).$
\item"(c)" Any two  
half-twists     are conjugate to each other.
It is possible to choose the conjugacy element  such that it commutes
with  any
half-twist which is disjoint from the given two half-twists. \item"(d)"
Any two pairs of
disjoint (or transversal, consecutive, cyclic) half-twists are conjugate
to each other. 
\item"(e)" Any two good quadrangles   are conjugate to each other.
 \item"(f)"
Every two pairs of orderly adjacent (non-orderly adjacent) consecutive
half-twists are
conjugate to each other,   preserving the polarization.\endroster
\ep \demo{Proof} 
If the statements hold for $B_n$, they also hold for $\tilde B_n.$ We
shall prove them for
$B_n.$ \roster\item"(a)" By simple geometric observation in $(D,k).$
\item"(b)" By Claim 1.0(i),(ii).
\item"(c)" We consider here only 2 cases: The half-twists $X,Y$ are 
consecutive; the
half-twists   $X,Y$ are disjoint.
For $X$ and $Y$ consecutive we get from (b) that $(X)_{YX}=Y.$
For $X=H(x)$ and $Y=H(y)$ disjoint half-twists, we choose a simple path
$\sigma$ connecting
an end point of $x$ with an end point of $y$ and do not meet any other
point of $K$ and any
other point of $x\cup y.$
Denote $Z=H(\sigma)$.
Since $X$ is consecutive to $Z$ and $Z$ is consecutive to $Y$ we get
from the first case that
$(X)_{ZXYZ}=(X_{ZX})_{YZ}=Z_{YZ}=Y.$
\item"(d)" We only consider here the case in which $(X_1,X_2)$ and
$(Y_1,Y_2)$ are 2 pairs of
disjoint half-twists.
By (c) there exist $b_1$ such that $(X_1)_{b_1}=Y_1$ and $b_2$ such that
$(X_2)_{b_2}=Y_2.$
Since $X_1$ is disjoint from $X_2$, for $b_1$ chosen as in (c) we get
$(X_2)_{b_1}=X_2.$
Similarly, $(Y_1)_{b_2}=Y_1.$
Let $b=b_1b_2.$
We get $(X_1)_b=Y_1$ and $(X_2)_b=Y_2.$
\item"(e)" Without loss of generality we can assume that the quadrangles
are disjoint (the
paths and the interior).
Take 2 disjoint quadrangles as in Fig. 1.3(a). 

\midspace{1.80in}\caption{Fig.   1.3($a$)\qquad\qquad\qquad\qquad\qquad
Fig. 1.3($b$)}

\flushpar Let $z_i$\ $(i=1,3)$ be a simple path connecting the end point
of $x_i$ with the
beginning point of $y_i$\ $(i=1,3),$ see Fig. 1.3(b).
Let $b_1=Z_1X_1Y_1Z_1,$\ $b_3=Z_3\1X_3\1Y_3\1Z_3\1.$
The choice of $b_1$ assures that $(X_1)_{b_1}=Y_1,$\ $(X_3)_{b_3}=Y_3$
(see (c)).
Moreover, $(X_3)_{b_1}=X_3,$\ $(Y_1)_{b_3}=Y_1.$
Let $b=b_1b_3.$
Clearly, $(X_1)_b=Y_1$ and $(X_3)_3=Y_3.$
We still have to show that $(X_2)_b=Y_2,$\ $(X_4)_b=Y_4.$
Now, $(X_2)_b=(X_2)_{Z_1X_1Y_1Z_1Z_3\1X_3\1Y_3\1Z_3\1}.$
By Claim 1.0(iii) $(X_2)_{Z_1}=H(t_1)$ where $t_1$ connects the 2 end
points of $X_2$ and
$Z_1$ that do not intersect (see Fig. 1.3(b)).
By Claim 1.0(i), $H(t_1)_{X_1}=H(t_1).$
By Claim 1.0(iii), $H(t_1)_{Y_1}=H(t_2).$
Alternately, we apply Claim 1.0(i) and Claim 1.0(iii) and finally get
$(X_2)_b=Y_2.$
Similarly, we get $(X_4)_b=Y_4.$
\item"(f)" The choices made in (b)   preserve the polarization .
 \qquad\qed\endroster  \edm

Lemma 1.2(b) is an important relation in $\tilde B_n;$ its discovery was
the first step in
understanding the structure of $\tilde B_n.$

\proclaim{Lemma   1.2}\ If $\{{Y}_i\}\ \ i=1,\dots, 4$ is
 a good quadrangle in $B_n,$ then

\rom{(a)}\ \ $\tilde{Y}_1
\tilde{Y}_3 =\tilde{Y}_3 \tilde{Y}_1,$

\rom{(b)}\ \ $Y_1^2Y_3^2Y_4^{-2}Y_2^{-2}$ normally generates
$\ker(B_n\to \tilde B_n).$
In particular,\linebreak $\tY_1^2\
\tY_3^2=\tY_2^2\ \tY^2_4.$
\ep

\demo{Proof}

(a) Since $\tilde Y_1$ and $\tilde Y_3$   are disjoint half-twists,
Claim 1.1(d)
implies that they commute.  

(b)  Since all pairs of transversal half-twists are conjugate to each
other,
$\ker(B_n\to\tB_n)$ is normally generated by any such commutator or a
conjugate of it.
So to \newline prove (b) it is enough to show that
$Y_1^2Y_3^2Y_4^{-2}Y_2^{-2}$ is conjugate to
a commutator of transversal half-twists.

Let $X_1, X_2, X_3$ be 3 half-twists such that $X_1$
and $X_2$ are consecutive, $X_2$ and $X_3$ are consecutive and $X_1$ and
$X_3$ are disjoint.
Denote $X_i = H(x_i),$ \linebreak $ i=1,\dots, 3.$
We use Claim 1.0(iii) to get a geometric description of \linebreak
$(X_3)_{X_2\1}=H((x_3)X_2\1)\ (=(X_2)_{X_3})$ and of
$(X_1)_{X_2\1}=H((x_1)X_2\1)\
(=(X_2)_{X_1}),$ (see Fig. 1.4(a)). We thus get that  $X_1,
(X_3)_{X_{2}\1}, X_3,
(X_1)_{X_{2}\1}$ is a good quadrangle.

 By Claim 1.0(iv),
$(X_3^{-2})_{X_2\1}(X_1)^{-2}_{X_2\1}X_1^2X_3^2=[X_2,(X_2)_{X_{1}X_{3}}]$.
By Claim 1.0(iii) applied twice,
  $(X_2)_{X_1X_3}$ is transversal to $X_2$
(see  Fig. 1.4(b) and (c)). Therefore
$(X_3)\1_{X_2\1}(X_1)^{-2}_{X_2\1}X_1^2X_3^2$ is a
commutator of transversal half-twists.

 \midspace{1.00in}\caption{Fig.   1.4($a$)}

\midspace{1.20in}\caption{Fig.   1.4($b$)\qquad\qquad\qquad\qquad Fig.  
1.4($c$)}

Since every 2 quadrangles are conjugate to each other,
$Y_2^{-2}Y_4^{-2}Y_1^2Y_3^2$ is
conjugate to a commutator of transversal half-twists.
As explained in the beginning this is enough to conclude that 
$Y_1^2Y_3^2Y_4^{-2}Y_2^{-2}$
normally generates $\ker (B_n\to \tB_n).$ In particular,
$\tY_1^2\tY_3^2\tY_4^{-2}\tY_2^{-2}=1$ and $\tY_1^2\ \tY_3^2=\tY_2^2\
\tY^2_4.$
 \qed \edm

\bigskip

\subheading{\S2. $\bold{\tB_n}$-groups and prime elements}

In this section we define a $\tilde B_n$-group, prime elements of
$\tilde B_n$-groups and
their supporting half-twist.  We present basic properties of a prime
element (Lemma 2.1) and
we prove 2 criteria for an element to be prime (Lemma 2.2 and Lemma
2.4).
These criteria are necessary steps for Proposition 7.1 in which we prove
a criterion for
prime elements which is easier to apply.
These elements are called prime elements since they satisfy an existence
and a uniqueness
property as proven in Theorem 3.3.

\definition{Definition}\ $\un{\tB_n-\text{group}.}$

A group $G$ is called a $\tB_n$-group if there exists a homomorphism
$\tB_n\ri\Aut (G).$
For $g\in G$ and $b\in \tB_n$ we denote by $g_b$ or $(g)_b$, the element
obtained from $b$
acting on $g.$  \enddefinition

\definition{Definition}\ $\un{\text{Prime element, supporting
half-twist, 
corresponding central element,}}$

\quad \quad\quad\quad$\un{\text{axioms of prime elements}}$

Let $G$ be a $\tB_n$-group.

An element $g\in G$ is called a prime element of $G$ if there exists a
half-twist $X\in B_n$ and $\tau\in \Center(G)$ with $\tau^2=1$ and
$\tau_b=\tau\ \forall\ b\in \tilde B_n$ such that:
\roster
\item $g_{\tX\1}=g\1 \tau.$
\item For every half-twist $Y$ adjacent to $X$ we have: 
\item"" (a)\ $g_{\tX\tY\1\tX\1}=g_{\tX}\1 g_{\tX\tY\1}$ 
\item"" (b)\ $g_{\tY\1
\tX\1}=g\1 g_{\tY\1}.$
\item For every half-twist $Z$ disjoint from $X,$\ $
g_{\tZ}=g.$\endroster

The half-twist $X$ (or $\tX$) is called the {\it{supporting half-twist}}
of $g,$

The element $\tau$ is called the {\it{corresponding central element}}.

We shall refer to (1), (2) and (3) as {\it Axiom} (1), {\it Axiom} (2),
{\it Axiom} (3),
respectively.\enddefinition 

The following lemma gives the basic properties of prime elements.
\proclaim{Lemma   2.1}

Let $G$ be a $\tB_n$-group.

Let $g$ be a prime element in $G$ with supporting half-twist $X$ and
corresponding central element $\tau.$
Then:
\roster
\item $g_{\tX}=g_{\tX\1}=g\1\tau,$\quad $g_{\tX^2}=g.$
\item $g_{\tY^{-2}}=g\tau \quad \forall \ \ Y$ consecutive half-twist to
$X.$
\item $[g, g_{\tY\1}]=\tau \quad \forall \ \ Y$ consecutive half-twist
to $X.$
\endroster\ep

\smallskip

\demo{Proof}
\roster
\item $g_{\tX^{-2}} = (g_{\tX\1})_{\tX\1} = (g\1\tau)_{\tX\1} =
(g\1)_{\tX\1}\cdot\tau = (g\1\tau)\1 \tau=g\Ri
g_{\tX}=g_{\tX\1}\overset\text{Axiom (1)}\to=g\1\tau.$
 \item $g_{\tX\1 \tY\1 \tX\1}=
(g_{\tX\1})_{\tY\1\tX\1}= (g\1\tau)_{\tY\1\tX\1}=
(g_{\tY\1\tX\1})\1\cdot\tau
\overset\text{Axiom (2)}\to =
g\1_{\tY\1}\cdot g\cdot \tau.$ \endroster
On the other hand,
$$\align g_{\tX\1 \tY\1 \tX\1}&=g_{\tY\1 \tX\1
\tY\1}=(g_{\tY\1\tX\1})_{\tY\1}\overset\text{Axiom (2)}\to =
 (g\1 g_{\tY\1})_{\tY\1}\\
&=  g\1_{\tY\1}\cdot g_{\tY^{-2}}.\endalign$$
Thus, $g_{\tY^{-2}}=g\tau.$
\roster
\item"(3)" $g_{\tX \tY\1 X\1} \overset{\text{Axiom (2)}} \to = 
g_{\tX}\1\cdot g_{\tX\tY\1} \overset{\text{by (1)}} \to = g\1_X
g_{\tX\tY\1}= g\cdot g\1_{\tY\1}\cdot \tau^2.$\endroster
On the other hand,
$$\align g_{\tX \tY\1 \tX\1}& \overset\text{by (1)}\to=
(g\1\tau)_{\tY\1\tX\1}\\ & \overset\text{Axiom (2)}\to = (g\1\cdot
g_{\tY\1})\1\cdot\tau = g\1_{\tY\1}\cdot g\cdot \tau.\endalign$$
Thus, $g\cdot g\1_{\tY\1}\cdot  = g\1_{\tY\1}\cdot g\tau.$
\flushpar Thus, $g_{\tY\1}\cdot g\cdot g_{\tY\1}\1\cdot g\1 =\tau\1 =
\tau.$
\flushpar Thus, $[g_{\tY\1}, g] = \tau.$ \qed \edm

\proclaim{Lemma   2.2}
Let $G$ be a $\tB_n$-group.
Let $g$ be a prime element in $G$ with supporting half-twist $X$ and
corresponding central
element $\tau.$
Let $b\in \tB_n.$
Then $g_b$ is a prime element with supporting half-twist $X_b$ and
central element $\tau.$
\ep
\pf We use the fact that $(a_b)_c=(a_c)_{b_c}$ and $(ab)_c=a_cb_c.$
We have to prove 3 properties:
\roster\item $g_{\tX\1}=g\1\tau \Ri (g_{\tX\1})_b=(g\1\tau)_b \Ri
g_{\tX\1b}=g_b\1\tau\Ri g_{bb\1\tX\1b}=(g_b)\1\tau\Ri
(g_b)_{\tX_b\1}=(g_b)\1\cdot \tau.$
\item Let $Y$ be a half-twist adjacent to $X_b.$
Then $Y_{b\1}$ is adjacent to $X$ and satisfies axiom (2) of prime
elements
for $g,$ $X$ and $Y_b\1.$
Namely: $g_{\tY_{b\1}\1\tX\1}=g\1g_{\tY_{b\1}\1}$ and
$g_{\tX\tY_{b\1}\1\tX\1}=g\1_{\tX}g_{\tX\tY_{b\1}\1}.$
 \item"" (a) \
$g_{\tY_{b\1}\1\tX\1}=g\1g_{\tY_{b\1}\1}\Ri$ \newline

$(g_{\tT_{b\1}\1\tX\1})_b=(g\1g_{\tY_{b\1}\1})_b\Ri$\newline

$(g_b)_{(\tY_{b\1}\tX\1)_b}=(g_b)\1(g_b)_{\tY\1}\Ri$\newline

$(g_b)_{\tY\tX_b\1}=g_b\1\cdot(g_b)_{\tY\1}.$ 
\item"" (b)\ $g_{\tX\tY_{b\1}\1\tX\1}=g\1_{\tX}g_{\tX\tY_{b\1}\1}\Ri$
\newline

$(g_b)_{\tX_{b}\tY\1\tX_b\1}=(g_b\1)_{\tX_b}(g_b)_{\tX_b\tY\1}.$
\item"(3)" Let
$Z$ be a half-twist disjoint from $X_b.$ Then $Z_{b\1}$ is disjoint from
$X.$
Then $g_{\tZ_{b\1}}=g.$
We conjugate $g_{\tZ_b\1}=g$ by $b$ to get:\
$(g_b)_{\tZ}=g_b.$\quad $\square$\endroster\enddemo

 We need the following technical lemma on $B_n$ to prove later a
criterion for a prime
element in a $\tB_n$-group. 
\proclaim{Claim   2.3}
Let $(X_1, \ldots, X_{n-1})$ be a frame in
$B_n = B_n[D,{K}]$. 
 Let $C(X_1) = \{b \in B_n|[b,X_1] =
1\}$ (centralizer of $X_1$), $C_p(X_1) = \{b \in
B_n\bigm|(X_{1})_{b} = X_1,$ preserving the polarization\}. 
Let
 $\sigma = X_2X^2_1X_2$.
Then $C(X_1)$ is generated by $\{X_1,\sigma,X_3, \ldots, X_n\}$,
$C_p(X_1)$ is generated by   $\{X^2_1, \sigma, X_3, \ldots,
X_n\}$.
\ep
\demo{Proof} Let $K=\{a_1,\dots,a_n\}.$  Let $x_1, \ldots ,
x_{n-1}$ be a system of
consecutive simple paths in $D,$ such that $X_i=H(x_i)$ ($H(x_i)$ is 
the half-twist corresponding to $x_i;\ 
x_i$ connects $a_i$ with
$a_{i+1}$). 
 Let $\Gamma_1, \ldots, \Gamma_n$ be a free geometric base
of $\pi_1(D-K,*)$ consistent with $(X_1, \ldots, X_{n-1})$ (that
is, $(\Gamma_{i+1})X_i = \Gamma_i$, $(\Gamma_i)X_i =
\Gamma_i\Gamma_{i+1}\Gamma^{-1}_i$, $(\Gamma_j)X_i = \Gamma_j$
for $j \neq i$, $i + 1$).  
We can assume that the $x_i$ do not
intersect the ``tails'' of $\Gamma_1, \ldots, \Gamma_n$.

Let $K_1$ be a finite set of $D$ obtained from $K\cup\{x_1\}$ by
contracting
$x_1$ to a point $\tilde a_2\in x_1.$
$K_1=\{\tilde a_2,a_3,\dots, a_n\}.$
Let $B_{n-1}=B_{n-1}[D,K_1].$
Let $Y_2,\dots ,Y_{n-1}$ be a frame of $B_{n-1},$ where $Y_i$ can be
identified with $X_i$ for $i\ge 3.$

Let $H=\{b\in B_{n-1}\bigm| (\tilde a_2)b=\tilde a_2\}.$
From the short exact sequence $1\rightarrow P_{n-1}\hookrightarrow
B_{n-1} \rightarrow S_{n-1} \rightarrow 1$ \ ($P_{n-1}$ the pure braid
group) we
can conclude that $H$ is generated by $Y_3,\dots ,Y_{n-1}$ and by the
generators of $P_{n-1}.$
We remove the generators of $P_{n-1}$ that can be expressed in terms of
$Y_3,\dots,Y_{n-1}$ (see \cite{A},\cite{B}, and \cite{MoTe1}, Section
IV) and
conclude that $H$ is generated by $Y_2^2,Y_3,\dots ,Y_n.$
The element $Y_2^2$ corresponds to the motion $\Cal M'$ of $\tilde
a_2,a_3,\dots,a_n $ described as follows:
\ $\tilde a_2,a_4,\dots, a_n$ stays in place and $a_3$ is moving around
$\tilde a_2$ in the positive direction (see Fig.   2.(a)).

\midspace{2.100in}\caption{Fig.   2}

We  define a homomorphism $\Phi: C_p(X_1)\ri H$ as
follows:

Let $U$ be a ``narrow'' neighborhood of $x_1$ such that $\lambda =
\partial U$ is a simple loop.  
Take $b \in C_p(X_1).$ 
It can be
represented by a diffeomorphism $\beta:  D \rightarrow D$  such that
$\beta({K}) = {K}$,
$\beta|_{\partial D} = \operatorname{Id}_{\partial D}$ and
$\beta|_{\overline{U}} = \operatorname{Id}_{\overline{U}}$
$(\overline{U} = U \cup \lambda)$. 

The diffeomorphism $\be$ also defines an element of $B_{n-1}[D,K_1].$
This element is in fact in $H$ since $\tilde a_2\in x_1$ and thus
$(\tilde a_2)\be=\tilde a_2.$
Denote this element by $\Phi(b).$
The map $\Phi$ constructed in this way is obviously a homomorphism,
$\Phi: C_p(x_1)\ri H.$
Clearly, $X_3,\dots, X_{n-1}\in C_p(X_1).$
Clearly, $\Phi(X_i)=Y_i$ for $i\ge 3.$
Let $\Cal M$ be the following motion   $(D,K)$: \
$a_1,a_2,a_4,\dots,a_n$ are stationary and $a_3$ goes around $a_1,a_2$
in the 
positive direction (Fig.   2(b)).
Let $u$ be the braid in $C_p(X_1)$ induced from the motion  $\Cal M.$
Clearly, $\Phi(u)=Y_2^2.$
Thus, $\Phi$ is onto and $\Phi(u),\Phi(X_3),\dots,\Phi(X_{n-1})$
generate
$H.$
One can check that $u=Z_{31}^2Z_{32}^2.$
But $Z_{31}=X_2X_1X_2\1$ and $Z_{32}=X_2.$
Thus $u=\si.$
Thus $C_p(X_1)$ is generated by $\si,X_3,\dots, X_{n-1}$ and a set of
generators for $\ker\Phi.$

 Consider $\pi_1(D-K\cup x_1,*).$ Let $\tilde
\G_2$ be the path obtained by connecting $\lambda$ with $*\in \partial
D$ by a simple path intersecting each of $\Gamma_3, \ldots, \Gamma_{n }$
only at $*$. We get a (free) geometric base $\tilde{\Gamma}_2,\Gamma_3,
\ldots, \Gamma_n$ of $\pi_1 (D - ({K} \cup x_1),*)$. 
It is obvious that
$\tilde\G_2=\G_1\G_2.$ \ $\Phi(b)$ defines in a natural way an
automorphism of $\pi_1(D_k\cup\{X_1\},*)$ such that
 $\Phi(b)$ does not
change the product $\tilde{\Gamma}_2 \Gamma_3 \ldots \Gamma_n$, and
$(\tilde{\Gamma}_2)\Phi(b)$ is a conjugate of $\tilde{\Gamma}_2$.

Consider now any $Z \in \ker\Phi$.  
We have
$(\tilde{\Gamma}_2) Z = \tilde{\Gamma}_2$ $(\tilde{\Gamma}_2 =
\Gamma_1\Gamma_2)$, $(\Gamma_j)Z = \Gamma_j \; \forall j = 3,
\ldots , n$.  
This implies that $Z$ can be represented by a
diffeomorphism which is the identity outside of $U$, that is, $Z =
X^\ell_1$,\ $\ell \in {\Bbb Z}$. 
 Since $Z \in C_p(X_1)$,\linebreak
  $\ell$ must equal $0\pmod 2$.

Thus, $C_p(X_1)$ is generated by $X^2_1,\sigma,
X_3, \ldots, X_{n-1}$.   
Clearly, $C(X_1)$ is generated by
$ {C}_p(X_1)$ and $X_1$.  \qed\quad for the Claim\edm

\proclaim{Lemma   2.4}
Let $\{X_1, \ldots, X_{n-1}\}$ be a frame in
$B_n$, $(\tX_1, \ldots, \tX_{n-1})$ their images in $\tilde{B}_n$. 
Let $G$ be a $\tB_n$-group.
Let $u \in G$,\ $\tau\in G$ be such that
\roster\item\  $u_{\tX^{-1}_1} = u^{-1} \tau$ with $\tau^2 = 1$, $\tau
\in$
$\Center(G)$, $\tau_b = \tau \; \forall b \in \tilde{B}_n$;
\item"(2$_a$)" \ $u_{\tX^{-1}_2\tX^{-1}_1} = u^{-1}u_{\tX^{-1}_2}$; 
\item"(2$_b$)" \
$u_{\tX_1\tX^{-1}_2\tX^{-1}_1} = u^{-1}_{\tX_1} u_{\tX_1\tX^{-1}_2}$;
\item"(3)"   $u_{\tX_j} = u \; \forall j = 3, \ldots, n - 1$.\endroster

Then $u$ is a prime element in $G$, and $X_1$ is its supporting
half-twist and $\tau$ is its central element.
\ep

\demo{Proof}
Let $Z \in B_n$ be any half-twist disjoint from $X_1$, $\tZ$ be the
image of $Z$ in $\tilde{B}_n$.  
$\exists b \in B_n$ such that
$(X_{1})_{b} = X_1$, $(X_{3})_{b} = Z$.  
By Claim    2.3, $b$  
belongs to the subgroup of $B_n$ generated by $X_1,X_3, \ldots, X_{n-1}$
and
$\sigma = X_2X^2_1X_2$. 
 Let $\tilde{b}$ and $\tilde{\sigma}$ be the images of
$b$ and $\sigma$ in $\tilde{B}_n$.   
We have $u_{\tilde\si\1}=
u_{\tX^{-1}_2\tX^{-2}_1\tX^{-1}_2} = (u_{\tX^{-1}_2\tX^{-1}_1})_{\tX^{-
1}_1\tX^{-1}_2} = (u^{-1}u_{\tX^{-1}_2})_{\tX^{-1}_1\tX^{-1}_2} = (\tau
u)_{\tX^{-1}_2} \cdot u_{\tX^{-1}_2\tX^{-1}_1\tX^{-1}_2} = \tau
u_{\tX^{-
1}_2} \cdot u_{\tX^{-1}_1\tX^{-1}_2\tX^{-1}_1} = \tau u_{\tX^{-1}_2}
\cdot (\tau u^{-1})_{\tX^{-1}_2\tX^{-1}_1} = \tau^2 u_{\tX^{-1}_2} u^{-
1}_{\tX^{-1}_2} \tau(\tau u) = u$.  
Then $u_{\tilde\si}=u.$ 
 Now: $u_{\tX_j}=u$
for $j\ge 3$ (by assumption (3)) and $u_{\tX_1^2}=u$ (by assumption
(1)). 
Thus, if $X_1$ appears in $b$ an even number of times, then $u_b=u.$
Otherwise
we replace $b$ by $bX_1.$
 The ``new'' $b$ satisfies the same requirement for
$b$ as above, as well as the equation $u_{\tilde b}=u.$ 
Thus, we can assume $u_{\tilde{b}} =
u$. 
 We have $$u_{\tZ} = u_{\tilde{b}^{-1}\tX_3\tilde{b}} =
u_{\tX_3\tilde{b}} =
u_{\tilde{b}} = u.$$

Let $Y$ be a half-twist in $B_n$ adjacent to $X_1$.  $\exists b_1
\in B_n$ such that $(X_{1})_{b_1} = X_1$, $(X_{2})_{b_1} = Y$. 
 Let
$\tilde{b}_1$ and $\tY$ be the images of $b_1$ and $Y$ in
$\tilde{B}_n$. 
 As above, we can choose $b_1$ so that
$u_{\tilde{b}_1} = u$.  
Applying $\tilde{b}_1$ on the assumptions (2$_a$) and
(2$_b)$  we get (using $u_{\tilde{b}_1} = u$,
$(\tX_{1})_{\tilde{b}_1} = \tX_1$,\ $(\tX_{2})_{\tb_1} = \tY$):
$$u_{\tY^{-1}\tX^{-1}_1} = u^{-1}u_{\tY^{-1}} \ \text{and} \
u_{\tX_1\tY^{-1}\tX^{-1}_1} = u^{-1}_{\tX_1} u_{\tX_1\tY^{-1}}.\quad
\square$$ \edm

\bk

\subheading{\S3. Polarized pairs and uniqueness of coherent pairs}

In this section we extend the notion of a polarized half-twist to prime
elements, and we
create polarized pairs which consist of a prime element and its
supporting half-twist when
considered polarized. We also extend the notion of braids preserving the
polarization to
coherent pairs.

The main results of this section are Theorem 3.3 and Corollary 3.5.
In Theorem 3.3 we establish the unique existence of a polarized pair,
with a given
supporting half-twist coherent to an original polarized pair.
We denote by $L_{h,\tilde X}(\tilde T)$ the unique prime element with
supporting half-twist
$\tilde T$ in the new polarized pair which is coherent to $(h,\tilde
X)$. In Corollary 3.5 we
prove simultaneous conjugation.

\demo{Definition}\ $\un{\text{Polarized pair}}$

Let $G$ be a $\tilde B_n$-group, $h$ a prime element of $G,$\ $X$
its supporting half-twist.
If $X$ is polarized, we say that $(h, X)$ (or $(h,\tX))$ is a polarized
pair
with central element $\tau,$\quad $\tau=hh_{\tX\1}.$\edm

\newpage

\demo{Definition}\ $\un{\text{Coherent pairs, Anti-coherent pairs}}$

We say that two polarized pairs $(h_1, \tX_1)$ and $(h_2, \tX_2)$ are
coherent
(anti-coherent) if $\exists \tb \in \tB_n$ such that
$(h_1)_{\tb}=h_2,$\ $(\tX_1)_b=\tX_2,$ and $b$ preserves (reverses) the
polarization.\edm 
\proclaim{Corollary   3.1} Coherent and anti-coherent polarized pairs
have
the same central element.\ep
\demo{Proof} The prime elements of coherent and anti-coherent pairs are
conjugate to each other.
Thus by Lemma   2.2 we get the corollary.\quad $\square$\edm

We need the following Lemma to prove later the unique existence of a
polarized pair with given
supporting half-twist coherent to an original polarized pair.

\proclaim{Lemma   3.2}
Let $G$ be a $\tilde B_n$-group.
Let $h\in G$ be a prime element with supporting half-twist $X.$
Let $b\in B_n.$
Then $X_{ b}= X$ preserving the polarization $\Rightarrow h_{\tilde
b}=h.$\ep

\demo{Proof}
  We can choose a set of
standard generators for $B_n[D,{K}],$ $\{X_1, \ldots,
X_{n-1}\}$ with $X_1 = X$. 
 Let $\sigma=X_2X_1^2X_2.$ 
Consider $C_p(X_1),$ the
centralizer of $X_1,$ preserving  the polarization.
 By Lemma   2.3, $C_p(X_1)$
is the subgroup of $B_n$ generated by $X^2_1, \sigma, X_3, \ldots,
X_{n-1}$.   
Since  $\tX_3,
\ldots, \tX_{n-1}$   are disjoint from $\tX_1$, they do  not change $h$
(by 
axiom(3) of prime elements).  
By Lemma   2.1, $h_{\tX_1^2}=h$.  
Consider
 $h_{\tilde\sigma^{-1}} = h_{\tX^{-1}_2\tX^{-2}_1\tX^{-
1}_2}$.  We have:
$$\align
h_{\tilde\sigma^{-1}} &=
h_{\tX^{-1}_2\tX^{-2}_1\tX^{-1}_2}= (h_{\tX^{-
1}_2\tX^{-1}_1})_{\tX^{-1}_1
\tX^{-1}_2}\overset\text{by Axiom(2) of prime
element}\to  =
(h^{-1}h_{\tX^{- 1}_2})_{\tX^{-1}_1\tX^{- 1}_2} \\
  & = 
(h^{-1}_{\tX^{-1}_1}h_{\tX^{-1}_2\tX^{-1}_1})_{\tX^{-1}_2}\overset\text{by
Axiom(2) 
of prime
element}\to =
(\tau h \cdot h^{-1}h_{\tX^{-1}_2})_{\tX^{-1}_2} = \tau h_{\tX^{-2}_2}\\
&\overset\text{by Lemma   2.1(2) }\to
  =  \tau h \tau = h.\endalign$$
Thus $h_{\tilde\sigma} = h.$
Thus, for every generator $g$ of $C_p(X),$\ $h_{\tg}=h.$ Since
$ b \in  {C}_p(X),$ \ $h_{\tilde b} = h$.\quad
$\square$\edm
\proclaim{Theorem  3.3} \rom{(Existence and Uniqueness Theorem)}\
Let $G$ be a $\tilde B_n$-group.
Let $(h,\tX)$ be a polarized pair of $G$.
 Let $\tT$ be a polarized half-twist in
$\tilde{B}_n$.  
Then there exists a {\it unique} prime element
$g \in G$ such that $(g,\tT)$ and
$(h,\tX)$ are coherent.
\ep

\demo{Proof} Let $X,T \in B_n$ be polarized half-twists
representing $\tX$ and $\tT$.  $\exists b \in B_n$ such that $T =
X_b$ preserving the polarization. 
 Let $\tilde{b}$ be the image of
$b$ in $\tilde{B}_n$.  
Taking $g = h_{\tilde{b}}$ we
obtain a polarized pair $(g,\tT)$ coherent with  $(h,\tX)$.  
To prove the {\it
uniqueness} of $g$, assume that
$(g_1,\tT)$ is another polarized pair coherent with $(h,\tX)$.  
Then
$\exists b_1 \in B_n$ with $g_1 = h_{\tilde{b}_1}$  and $T = X_{b_1}$,
preserving the polarization. 
 We have $T = X_{b_1} = X_b$ and
$X_{b_1b^{-1}} = X$. 
 Denote  $b_2 = b_1b^{-1}$, so $X_{b_2} =
X$ (preserving the  polarization). 
By the previous lemma, $h_{\tilde b_2}=h.$
Thus, $h_{\tilde b_1}=h_{\tilde b}$ or $g=g_1.$\quad $\square$\edm
\demo{Definition} $\un{L_{(h,\tX)}(\tT)}$

Let $(h,\tX)$ be a polarized pair.
Let $\tT\in\tB_n.$
We denote by $L_{(h,\tX)}(\tT)$   the unique prime element such that
$(L_{(h,\tX)}(\tT),\tT)$
is coherent with $(h,\tX).$\edm
\mk
\proclaim{Lemma   3.4}
Assume $(h,\tX) $ and $(g,\tX)$ are polarized pairs.  Let $\tau$ be the 
central
element of $(g,\tX).$ If $(h,\tX)$ is anti-coherent to $(g,\tX)$ then
$h=g\1\cdot\tau.$\ep

\demo{Proof}
By assumption, $\exists b\in B_n$ such that $g=h_{\tilde b}$ and
$X=X_b,$ reversing polarization. Thus $X_{bX\1}=X,$ preserving the
polarization.
Thus $(h_{\tilde b\tX\1},\tX)$ is coherent with $(h_{\tilde
b\tX\1},\tX_{\tilde b\tX\1}).$
Clearly, $(h,\tX)$ is coherent with $(h_{\tilde
b\tX\1},\tX_{\tilde b\tX\1}).$
From uniqueness, $h=h_{\tilde b\tX\1}=g_{\tX\1}.$
Since $\tau$ is the central element of $(g,\tX),$\ $\tau=g\cdot
g_{\tX\1}.$
Thus, $g_{X\1}=g\1\tau.$
So $h=g\1\tau.$\quad\qed\edm
From uniqueness we get  simultaneous conjugation:

 \proclaim{Corollary   3.5} If
$(a_i,\tX)$ is coherent with $(g_i,\tY)$\ $i=1,2$, then there exist
$b\in\tB$
such that $(a_i)_b=g_i,$\ $i=1,2.$\ep
 \demo{Proof} Let $b$ be the element of $\tB_n$
such that $(a_1)_{b}=g_1,$\ $(\tX)_{b}  =\tY.$ Now,
$((a_2)_{b},(X)_{b})$ is coherent
with $(a_2,\tX).$ Since $(\tX)_b=\tY,$\ $((a_2)_b,\tY)$ is coherent with
$(a_2,\tX).$ The pair $(g_2,\tY)$ is also coherent with $(a_2,\tX).$
From
uniqueness, $(a_2)_{b}=g_2.$ \qed\edm \bk

\subheading{\S4. $\bold{\tB_n}$-action of  half-twists}

In this section we study the action of half-twists on prime elements.
We compute it for the case when the half-twist is adjacent to the
supporting half-twist
(Proposition 4.1) and for the case in which it is transversal (Lemma
4.2).
The action of a disjoint half-twist is part of the defining axioms of
prime elements
(Axiom~(iii)).

\proclaim{Proposition   4.1} Let $G$ be a $\tilde B_n$-group.
Let $(h,\tX)$ be a polarized pair of $G$ with coresponding central
element $\tau.$
Let $T,Y$ be $2$ orderly adjacent polarized half-twists in
$ {B}_n$.  
 Denote by $\tY^\prime$ the polarized half-twist
obtained from $\tY$ by changing polarization (that is, $\tT,
\tY^\prime$ are {\it  not orderly adjacent}). 
 Denote by
$$\gather L(T) = L_{(h,\tX)}(\tT), \\
L(Y) = L_{(h,\tX)}(\tY),\\
L(Y^\prime) = L_{(h,\tX)}(\tY^\prime).\endgather$$
Then
\roster\item $L(T)_{\tT^{-1}} = L(T)^{-1}\tau;$ 
\item  $L(T)_{\tY^{-1}} = L(T)L(Y)$;
\item  $L(T)_{(\tY'){}^{-1}} = L(Y')^{-1}L(T)$.\endroster
\ep

\demo{Proof}

(1) By Lemma   2.1(1).

(2)  Let $b\in B_n$ be such that $L(T)=h_{\tilde b}$ and  $T=X_b,$
preserving the polarization.
Let $Y_1 = Y_{{b}^{-1}}.$  Then $X,Y_1$  are
adjacent half-twists $(X=T_{b\1},\ Y_1=Y_{b\1}),$ and so
$h_{\tY^{-1}_1\tX^{-1}} =
h^{-1}h_{\tY^{-1}_1}$.  Applying  $\tilde{b}$ to that equation, we get
$(L(T))_{\tY^{-1}\tT^{-1}} = L(T)^{-1} L(T)_{\tY^{-1}}$, or $$
L(T)_{\tY^{-1}} = L(T) \cdot L(T)_{\tY^{-1}\tT^{-1}}
$$

Let $b_1 = bY^{-1}T^{-1}$.  Then $X_{b_1} = X_{bY^{-1}T^{-1}} =
T_{Y^{-1}T^{-1}} = Y$ ($T_{Y\1}=Y_T$ since $Y,$ $T$ are adjacent). 
Since $T,Y$ are
orderly adjacent and $X_b = T,$ preserving the polarization, one can
easily check
that actually $X_{b_1} = Y,$ preserving the polarization. 
Due to {\it
uniqueness} of $L_{(h,\tX)}(\tY)=L(Y)$, we get
$L(Y) = h_{\tilde{b}_1} = h_{\tilde{b}\tY^{-1}\tT^{-1}} =
L(T)_{\tY^{-1}\tT^{-1}}$.  Together with the previous equation we get:
$L(T)_{\tY^{-1}} = L(T)L(Y)$, which is (2).

(3)  From (2) we get
$  L(T)_{\tY'{}^{-1}} = L(T)_{\tilde Y\1}=L(T)L(Y)$.
By Lemma 2.1 $[L(T),L(T)_{\tY^{-1}}] = \tau.$ 
Using  (2) we substitute $L(T)L(Y)$ instead of
$L(T)_{\tY\1}$   and get  $\tau = [L(T),L(T)L(Y)] = [L(T),L(Y)]$  and
$L(T)L(Y) = \tau L(Y)
L(T)$. Thus $L(T)_{\tY'{}\1}=\tau L(Y)L(T).$ 
 Using $Y_{Y^{-1}} = Y'$ (preserving the polarization) and 
uniqueness, we can write
$L(Y') = L(Y)_{\tY^{-1}}.$  By (1), $L(Y)_{\tY\1} = L(Y)^{-1}
\tau$, and thus $L(Y')=L(Y)\1\tau.$
As a consequence, we get $L(T)_{\tY'{}\1}=L(Y')\1L(T),$
 which is (3). \qed  \edm

\proclaim{Lemma   4.2} Let $G$ be a $\tilde B_n$-group.
Let $h\in G$ be a prime element with
supporting half-twist $X.$ 
Let $Z$ be a half-twist in $ {B}_n$
transversal to $X$.  Then $h_{\tZ} = h$.\ep

\demo{Proof}   Let $x,z$ be 2 transversally intersecting 
simple paths corresponding to $X,Z$ (see Fig.   4).

\midspace{2.00in}
\caption{Fig.   4}

There exists a simple path $y$ such that the corresponding half-twist
$Y$ is
adjacent to $X$ and $Z$, and $Z_1 = Z_{Y^{-2}}$ is disjoint from $X.$
Let $z_1$
be the path corresponding to $Z_1$ (see Fig.   4). 
 We have $h_{\tZ} =
h_{\tY^{-2}\tZ_1\tY^2} \overset\text{by Lemma   2.1}\to= (h
\tau)_{\tZ_1\tY^2} = (h \tau)_{\tY^2} = h \tau \cdot \tau = h$. 
\qed\quad\edm

\newpage

\subheading{\S5. Commutativity properties}

In this section we compute the commutator of prime elements depending on
the relative location
of their supporting half-twists.
We also prove a commutativity result for $\tilde P_n$ (the image of
$P_n$, the pure braid
group, in $\tilde B_n)$ to be used in studying the structure of $\tilde
P_n$ as a $\tilde
B_n$-group in Section 6 (see Lemma 6.1).

\proclaim{Proposition   5.1}  Let $G$ be a $\tilde B_n$-group.
Let $(g_1,\tY_1),$ \ $(g_2,\tY_2)$ be $2$ polarized pairs of $G$. 
Assume that they are coherent or anti-coherent.
Let $\tau$ be the corresponding central element of $(g_1,\tY_1)$
\ $(\tau = g_1(g_{1})_{\tY^{-1}_1})$.  

Then

\roster\item  if $\tY_1,\tY_2$ are adjacent, then $[g_1,g_2] = \tau$;
\item  if $\tY_1,\tY_2$ are disjoint or transversal, then $[g_1,g_2] =
1$.\endroster
\ep

\demo{Proof} We use in the proof a few identities concerning commutators
in a group:

$[a,bc]=[a,b][a,c]_{b\1}$

$[a,ba]=[a,b\1]_{b\1}\1=[a,b]$

$[a,a\1b]=[a,b]_a=[b,a\1]$

(1)  Assume first that $(g_1,\tY_1),(g_2,\tY_2)$ are coherent. 
Take $b \in \tilde{B}_n$ with $g_2 = (g_1)_{b}$, $\tY_2 = (\tY_1)_{b}$
(preserving the polarization).  Let $b_1 = \tY^{-1}_2\tY^{-1}_1$.  By
Claim 1.1(b)
$(\tY_1)_{b_1} = \tY_2$.  Assume that $b_1$ preserves the polarization
of
$\tY_1,\tY_2$.  We have $((g_{1})_{b_1},\tY_2)$ and $(g_2,\tY_2)$
coherent with $(g_1,\tY_1)$.  By Theorem   3.3 (the {\it
uniqueness} part) we get $(g_1)_{b_1} = g_2$.  Thus we have $g_2 =
(g_1)_{b_1} = (g_1)_{\tY^{-1}_2\tY^{-1}_1} \overset \text{by Axiom
(2)}\to= g^{-1}_1
(g_{1})_{\tY^{-1}_2}$, and $[g_1,g_2] =
[g_1,g^{-1}_1(g_{1})_{\tY^{-1}_2}] =  [g_1,(g_{1})_{\tY^{-1}_2}]_{g_1}
\overset\text{by Lemma   2.1(3)}\to=  \tau_{g_1} = \tau$.

If $b_1$ does not preserve the polarization of $\tY_1,\tY_2,$ then $b_2
=
b_1\tY_2$ does preserve it.    As
above, we get $g_2 = (g_1)_{b_2} = (g_1)_{\tY
^{-1}_2 \tY^{-1}_1\tY_2} = (g_1)_{\tY_1\tY^{-
1}_2\tY^{-1}_1} = (g^{-1}_1 \tau)_{\tY^{-1}_2\tY^{-1}_1} = \tau (g^{-
1}_{1})_{\tY^{-1}_2}g_1$, and then $[g_1,g_2] = [g_1,\tau (g^{-
1}_{1})_{\tY^{-1}_2}g_1] = [g_1,(g^{-
1}_{1})_{\tY^{-1}_2}] = [g_1,(g_1)_{\tY^{-
1}_2}]^{-1}_{(g_1)_{\tY^{-1}_2}}$\linebreak $ \overset\text{by Lemma  
2.1(3)}\to= \tau$.

Assume now that $(g_1,\tY_1)$, $(g_2,\tY_2)$ are anti-coherent.
Denote by $\tY'_2$
the half-twist obtained from $\tY_2$ by changing polarization.  One
can then check that $(g_1,\tY_1)$, $((g_{2})_{\tY^{-1}_2},\tY'_2)$ are
coherent.  By the above  $\tau = [g_1,(g_{2})_{\tY^{-1}_2}].$ 
By Corollary   3.1, $\tau$ is also the central element of $(g_2,\tY_2).$
Thus $\tau=g_2(g_2)_{\tY_2\1},$ which implies
$(g_2)_{\tY_2\1}=g_2\1\tau.$
Thus $\tau=[g_1,(g_2)_{\tY_2\1}] =
[g_1,g^{-1}_2 \tau] = [g_1,g_2^{-1}] = [g_1,g_2]^{-1}_{g_2}$.  Thus,
$[g_1,g_2]
= \tau^{-1}_{g_2} = \tau$. 

(2)  We can assume that $(g_1,\tY_1),$\ $(g_2,\tY_2)$ are coherent. 
(Otherwise, we replace $\tY_2$ by $\tY'_2$ and $g_2$ by $(g_{2})_{\tY^{-
1}_2}$ and use $[g_1,(g_{2})_{\tY^{-1}_2}] = [g_1,g_2]^{-1}_{g_2}$ as in
(1).)

Consider first the case where $\tY_1,\tY_2$ are disjoint.  We can
choose a standard base of $\tilde{B}_n$, say $(\tX_1,\tX_2, \ldots,
\tX_{n-1})$ such that $\tX_1 = \tY_1$, $\tX_3 = \tY_2$ and the given
polarizations of $\tY_1,\tY_2$ coincide with ``consecutive''
polarizations of $\tX_1,\tX_3$ (``end'' of $\tX_1$ = ``origin'' of
$\tX_2$, ``end'' of $\tX_2 =$ ``origin'' of $\tX_3$).  Let $b_1 = \tX^{-
1}_2\tX^{-1}_1\tX^{-1}_3\tX^{-1}_2$.  
By the proof of Claim 1.1, $\tY_2 = (\tY_1)_{b_1}$, preserving the
polarization. 
Then $((g_1)_{b_1},Y_2)=((g_1)_{b_1},(Y_1)_{b_1})$ is coherent with
$(g_1,Y_1).$
 From Theorem   3.3 ({\it uniqueness}) it follows that
$g_2 = (g_1)_{b_1} = (g_1)_{\tX^{-1}_2\tX^{-1}_1\tX^{-1}_3\tX^{-1}_2} =
(g^{-1}_1(g_1)_{\tX^{-1}_2})_{\tX^{-1}_3\tX^{-1}_2}$ \linebreak $ =
(g^{-1}_{1})_{\tX^{-
1}_2}(g_1)_{\tX^{-1}_2\tX^{-1}_3\tX^{-1}_2} = (g^{-1}_{1})_{\tX^{-
1}_2}(g_1)_{\tX^{-1}_3\tX^{-1}_2\tX^{-1}_3} = (g^{-1}_{1})_{\tX^{-
1}_2}(g_1)_{\tX^{-1}_2\tX^{-1}_3}$.    We can write $[g_1,g_2]
= [g_1,(g^{-1}_{1})_{\tX^{-1}_2}(g_1)_{\tX^{-1}_2\tX^{-1}_3}] =
[g_1,(g^{-1}_{1})_{\tX^{-1}_2}] \cdot [g_1,(g_1)_{\tX^{-1}_2\tX^{-
1}_3}]_{(g_1)_{\tX^{-1}_2}}.$
By   2.1(3)
$[g_1,(g_{1})_{\tX^{-1}_2}] = \tau$, which implies
$[g_1,(g_{1})_{\tX^{-1}_2\tX^{-1}_3}] = [g_1,(g_{1})_{\tX^{-
1}_2}]_{\tX^{-1}_3} = \tau_{\tX^{-1}_3}  = \tau$.
We also have
$[g_1,(g_1\1)_{X_2\1}]=[g_1,(g_1)_{X_2\1}]\1_{(g_1)_{X_2\1}}=\tau\1_{(g_1)_{X_2}\1}=\tau\1.$
Thus $[g_1,g_2]=\tau\cdot\tau\1=1.$

Assume now that $\tY_1,\tY_2$ are transversal.  As in the proof of
Lemma   4.2, we can find a half-twist $\tT \in \tilde{B}_n$ such that
$\tT$ is adjacent to $\tY_1,\tY_2$ and $ \tY_2' = (\tY_{2})_{\tT^{-2}}$
is
disjoint from $\tY_1$.  Let $b \in \tilde{B}_n$ be such that $\tY_2 =
(\tY_1)_{b}$, $g_2 = (g_1)_{b}$. 
Let
 ${g}_2' = (g_1)_{bT^2} = (g_{2})_{\tT^{-2}}$.  Then
$(g_2',\tY_2')$ is coherent  with
$(g_1,\tY_1)$.  Since $\tilde{Y}_2',\tY_1$ are disjoint, we get from
 the above $[g_1,g_2'] = 1$.  By Lemma   2.1 $ g_2' =
(g_2)_{\tT^{-2}} = g_2\tau$, or $g_2 = g_2'\tau$.  Therefore,
$[g_1,g_2] = [g_1, g_2' \tau] = [g_1,  g_2'] = 1$. \quad
\qed\edm

\demo{Note} Recall that there exists a natural homomorphism $\psi_n:
B_n\ri S_n$ where
$\psi_n(X_i)$ is the transposition $(i\ i+1)$ for $\{X_i\}_{i=1}^{n-1}$
a frame of $\tilde
B_n.$ 
Its kernel $P_n=\ker\psi_n$ is called the pure braid group.
Recall from \cite{B} and  \cite{ MoTe1} that  $P_n$ is generated by
$\{Z_{ij}^2\},$
where
$Z_{ij}=(X_i)_{X_{i+1}\dots X_{j-1}}.$\edm

\demo{Definition}\ $\un{\tilde{P}_n}$ 

$ \tilde{P}_n =
\ker(\tilde{B}_n\overset \tilde\psi_n\to \rightarrow S_n),$ where
$\psi_n$ is
induced naturally from $\psi_n\ (n\ge 4).$\edm

\proclaim{Proposition   5.2}
Assume $n \geq 4$. Let $\tX_1,\tX_2$ be $2$ adjacent half-twists in
$\tB_n.$
Let $c=[\tX_1^2,\tX_2^2].$  Then the commutant $\tilde{P}^\prime_n$ of
$\tilde{P}_n$ is generated by $c$ where $c_b = c \;
\forall b \in \tilde{B}_n$, and $c^2=1.$
Moreover, if  $(\tZ_1,\tZ_2)$ is another   pair  of adjacent
half-twists, then 
$[\tZ_1^2,\tZ_2^2]=[\tZ_1^2,\tZ_2\2]=[\tZ_1\2,\tZ_2\2]=c.$
In particular, if $(\tY_1,\tY_2)$ and $(\tZ_1,\tZ_2)$ are two pairs of
adjacent half-twists
then $[\tZ_1^2,\tZ_2^2] =[\tY_1^2,\tY_2^2].$ \ep

\pf Let $B_n = B_n[D, {K}]$). Complete $\tX_1$ and $\tX_2$ to 
 $\tX_1, \ldots, \tX_{n-1},$ a standard base of $\tilde{B}_n$, i.e.,
 $X_i=H(x_i)$ and $x_1,\dots,x_{n-1}$ are simple consecutive paths 
  in $D.$  Let $c=[\tX_1^2,\tX_2^2].$
 Let $X = (X_{1})_{\tX_2\tX_3}$.
By Claim 1.0(iii) repeated twice, $H=H(x)$ where $x$ is as in Fig.~5(a).
We have a quadrangle formed by $x_1,x_2,x_3,x$, (see
Fig.   5(a)).

\midspace{2.00in}
\caption{Fig.   5}

Evidently, $\tX_1,\tX_2,\tX_3,\tX$ form a good quadrangle in $\tB_n.$ 
Thus
by Lemma   1.2
$$
\tX^2_1\tX^2_3 = \tX^2_2\tX^2.
\tag1.10$$
Denote  $y_1 = \tX^2_1$,\quad $y_2 = \tX^2_2$,\quad $y_3 =
\tX^2_3$,\quad  $y_4
= \tX^2,$\ (the squares of the edges), $d_1 = \tX_1\tX^2_2\tX^{-1}_1$,\
$d_2 = \tX_2\tX_3^2\tX^{-1}_2$,\  (the squares of
the diagonals),  $y'= \tX_2\tX_1^2\tX^{-1}_2$ (the square of the outer
diagonal).
 See Fig.   5(b) where we denote the paths corresponding to the
half-twists whose squares we considered here.
Clearly,

$y_1y_3=y_3y_1$

$d_1=(y')_{y_1\1}$

$(y_3)_{x_2}=(d_2)_{x_3\2}=(d_2)_{y_3\1}$

$y_1y_2d_1=y_2y_1y'=\Delta_3^2$ (a central element of $P_3).$

We rewrite (1.10) to get
$$y_4=y_1y_3y_2\1
\tag1.11$$

Conjugating (1.11) by $\tX_1$, we get
$$
d_2 = y_1y_3y'{}^{-1}; 
\tag1.12$$
conjugating (1.11) by $\tX_2$ we get 
$y_4 = d_1(d_2)_{y_3\1}\cdot y_2\1.$ 
Since $d_1 =
y_1y'y^{-1}_1$, 

$$y_4 = y_1 y'y^{-1}_1 y_3d_2y^{-
1}_3y^{-1}_2 \overset\text{by (1.12)}\to=  y_1y' y^{-1}_1 y_3 y_1
y_3 y'{}^{-1} y_3^{-1}y^{-1}_2 .$$  
We compare the last expression with (1.11) to get $$y'y^2_3
y'{}^{-1}y^{-1}_3 =
y_3, \; \quad\text{or}\quad  [y',y^2_3] = 1.
\tag1.13$$

Since $y',y_3$ are squares of two adjacent half-twists in
$\tilde{B}_n$, and any two pairs of adjacent half-twists are
conjugate, we conclude from (1.13) that: 
$$ \text{If\ } \tilde Z_1,\tilde Z_2\ \text{is a pair of
adjacent half-twists in}\ \tilde{B}_n,\
[\tZ^2_1,\tZ^4_2] = 1,
\tag1.14$$
which also implies that $[\tZ^2_1,\tZ^2_2] = [\tZ^{-2}_1,\tZ^2_2] =
[\tZ^2_1,\tZ^{-2}_2] = [\tZ^{-2}_1,\tZ^{-2}_2].$

Conjugating (1.11) by $\tX^{-1}_3$ we get
$$d_1=y_1y_3(d_2)_{y_3\1} = y_1y_3\cdot y_3d^{-1}_2y^{-1}_3 = y_1y^2_3
\cdot
d^{- 1}_2y^{-1}_3\overset (1.14)\to = y_1d^{-1}_2y_3,$$

We have by (1.14) that $1 = [y^2_1,y_2] = [y_1,y_2]_{y^{-1}_1} \cdot
[y_1,y_2]$.  Denoting $c = [y_1,y_2]$, we can write
$$
c_{y_1^{-1}} = c^{-1},\quad\text{or}\quad c_{y_1} = c^{-1}.
\tag 1.15$$
  Denote by $\tilde{P}_3$ the
subgroup of $\tilde{B}_n$ generated by $y_1,y_2,d_1$, and by
$\alpha = \Delta^2_3 = y_1y_2d_1 = y_2y_1y'$ (a central
element of $\tilde{P}_3$), so that
$$
y' = y_1^{-1}y^{-1}_2 \alpha.
\tag 1.16$$
So, $c_{\tX_1} = [y_1,y_2]_{\tX_1} = [y_1,y']
\overset\text{by (1.16)}\to=  [y_1,y^{-1}_1y^{-1}_2 \alpha] = y_1 \cdot
y^{-1}_1y^{-1}_2 \alpha \cdot y^{-1}_1 \cdot \alpha^{-1}y_2y_1 =
y^{-1}_2y^{-1}_1y_2y_1 = [y^{-1}_2,y^{-1}_1] = [y_2,y_1] = c^{-
1}.$
Thus we have $c_{\tX^2_1} = (c^{-1})_{\tX_1} = c.$

 By (1.15)
$c_{\tX^2_1} = c_{y_1} = c^{-1}$.  

We compare the last two results to get\ $c = c^{-1}$ or
$$
c^2 = 1 \quad\text{and} \quad c_{\tX_1} = c.
\tag1.17$$

Using a conjugation which sends $(\tX_1,\tX_2)$ to $(\tX_2,\tX_1)$, we
obtain from (1.17)
$$
c^{-1}_{\tX_2} = c^{-1},\quad{or} \quad c_{\tX_2} = c.
\tag 1.18$$

(1.17) and (1.18) show that $\forall$  $z \in \tilde{B}_3$ (the
subgroup of $\tilde{B}_n$ generated by $\tX_1,\tX_2$) we have
$$
c_z = c.
\tag1.19$$

Consider now $c_{\tX_3} = [y_1,y_2]_{\tX_3} = [y_1,d_2]\overset\text{by
(1.12)}\to =
[y_1,y_1y_3y'{}^{-1}] = $\newline
$y_1y_1y_3y'{}^{-1}y^{-1}_1y' y^{-1}_3y^{-1}_1 = y_1y_3 \cdot
(y_1y'{}^{-1}y^{-1}_1 y')y^{-1}_3y_1^{-1}\overset\text{by (1.19)}\to= 
y_1 y_3 c y^{-1}_3y^{-1}_1 = c_{y^{-1}_1y^{-1}_3} = c_{y^{-3}_1} =
c_{y_3} =
c_{\tX^2_3}$, in short $c_{\tX^2_3} = c_{\tX_3}$.  This implies
$c_{\tX_3} = c$.

Since $c = [\tX^2_1,\tX^2_2]$, we have $\forall \tX_j$, $j \geq 4$,
$c_{\tX_j} = c$.  Therefore, $\forall j,$\ $c_{\tX_j}=c,$ and thus
$\forall b \in
\tilde{B}_n$ $c_b = c.$ 

Let $(\tY_1,\tY_2)$ be another pair of adjacent half-twists.
Since every 2 pairs of adjacent half-twists are conjugate in
$\tB_n,$\quad
$\exists b\in \tB_n$ such that
$[\tY_1^2,\tY_2^2]=[\tX_1,\tX_2^2]_b=c_b.$
Since\  $c_b=c\quad\fa b\in\tB_n,$\quad  $[\tY_1^2,\tY_2^2]=c.$
Since $c^2=1,$\ $c\in\Center(\tB_n),$ we also have
$[\tY_1^2,\tY_2\2]=[\tY_1\2,\tY_2\2]=c.$  In particular, if
$(\tZ_1,\tZ_2)$ is
another pair of adjacent half-twists, $[\tZ_1
^2,\tZ_2^2]=[\tY_1^2,\tY_2^2]=c.$
Because any two disjoint and transversal half-twists of $\tilde{B}_n$
commute,
and $\tilde{P}_n$ is generated by $\tZ^2_{ij} = (\tX^2_{i})_{\tX_{i+1}
\cdots
\tX_{j-1}},$ $1 \leq i < j \leq n$, we conclude that
$\tilde{P}^\prime_n$ is generated by $c$.  \quad \qed\enddemo

\remark{Remark 5.3} Throughout the rest of the paper we will use $c$ for
the generator of
$\tilde P_{n}'.$
We also recall from the above proposition  that $c=[\tX_1^2, X_2^2]$
for  any 2 consecutive
half-twists in $B_n$,\ $c\in \Center B_n$ and $c^2=1.$\endremark

\bigskip

\subheading{\S6. On the structure of $\bold{\tilde B_n}$ and
$\bold{\tilde P_n}$}

In Proposition 6.4   we establish a $\tilde
B_n$-isomorphism between $\tilde P_n$ and the $\tilde B_n$-group,
$\underline G(n),$ whose
definition and structure are given before the proposition (Claim 6.2).
Proposition 6.4 will be used in the final section when we give the
simplest criterion for a
prime element.
Corollary 6.5 establishes the fact that $\tilde B_n$ is an extension of
a solvable group by a
symmetric group, which is a structure type theorem for $\tilde B_n.$
 
If $G$ is a group, we denote its abelianization by $Ab(G).$
Recall that $\Ab(B_n)\simeq \Bbb Z$.
(This follows from the fact that  $B_n$ is generated by the half-twists 
and
every 2 half-twists are conjugate to each other (Claim 1.1(c)).
 \demo{Definition}\ $\un{P_{n,0}}$

 $P_{n,0} = \ker(P_n \rightarrow \Ab B_n)$ (``degree zero'' pure
braids).
\edm

\demo{Definition}\ $\un{\tilde{P}_{n,0}}$

$\tilde{P}_{n,0}$ is
the image of $P_{n,0}$ in $\tilde{P}_n$.  \edm

\remark{Remark 6.0}
Let $X_1,\dots,X_{n-1}$ be a frame of $B_n,$\
$Z_{ij}=(X_i)_{X_{i+1}\cdot\dots\cdot
X_{j-1}}.$
Since $  P_n$ is generated by $\{Z_{ij}^2\}_{i<j}$ (see \cite{B}),
$\tilde P_{n,0}$ is
generated by $\{\tX_1^2\tZ_{ij}^{-2}\}_{i<j}.$
If $\tilde g_1\dots\tilde g_m$ is another set of generators for $\tilde
P_{n,0},$ then
$\{\tX_1^2,\tilde g_1,\dots ,\tilde g_m\}$ generates $\tilde
P_n.$\endremark
 \proclaim{Lemma  
6.1} Let $X_1,X_2$ be $2$ consecutive half-twists in $B_n.$ Let\newline
$u =
(\tX^2_{1})_{\tX^{- 1}_2} \tX^{-2}_2$. Then $u\in\tP_{n,0}$,  \ $u$ is a
prime element in
$\tilde{P}_n$ (considered as a $\tilde{B}_n$-group),  $\tX_1$ is its
supporting half-twist
and $c=[\tX_1^2,\tX_2^2]$ is its coresponding central element.\ep

\demo{Proof} Clearly, $u \in \tilde{P}_{n,0}$.
We often use here the fact that $(X_1)_{X_2\1}=(X_2)_{X_1}$ as well as
the fact
that $[\tX_1^{\pm 2},\tX_2^{\pm 2}]=c,$ i.e.,
$\tX_1\2\tX_2^2=c\tX_2^2\tX_1\2$
and $(\tX_1^2)_{X_2\1}=(\tX_1^2)_{X_2}c$\quad for $c\in\Center(\tB_n),$\
$c^2=1.$
 In particular,
$u=(\tX_1)_{\tX_2\1}^2\tX_2\2=(\tX_2)_{\tX_1}^2\tX_2\2.$
Complete $X_1,X_2$ to a frame of $B_n:$\ $X_1,\dots,X_{n-1}.$ $(\la
X_i,X_{i+1}\ra=1$ and $[X_i,X_j]=1$\ $|i-j|>2.).$ We shall use Lemma  
2.4,
and so we must check  conditions (1), $(2_a)$, $(2_b)$, (3) of the
Lemma.

(1)  We have $u_{\tX^{-1}_1} = (\tX^2_{2})_{\tX_1\tX^{-1}_1} \cdot 
(\tX^{-
2}_{2})_{\tX^{-1}_1} =\tX^{2}_{2}\cdot (\tX_1\2)_{\tX_2}.$
Let $c = [\tX^2_1,\tX^2_2].$
By Proposition   5.2, $(\tX_1\2)_{\tX_2}=(\tX_1\2)_{\tX_2\1} c.$ Thus
$u_{\tX_1\1}=  \tX^2_2 \cdot (\tX^{-2}_{1})_{\tX^{-1}_2} c = u^{-1}c$.
By Proposition 5.2, we also know that $c_b=c$\ $\forall b\in \tB_n,$\
$c\in \Center \tB_n$
and $c^2=1.$

$(2_a)$ Since $[\tX_1^2,\tX_2^2]=c,$\ $u_{\tX^{-1}_2} =
(\tX^2_{1})_{\tX^{-2}_2}
\cdot \tX^{-2}_2 = c\tX^2_1\tX^{-2}_2 = \tX^{-2}_2 \cdot \tX^2_1$, and
$u_{\tX^{-1}_2\tX^{-1}_1} = (\tX^{-2}_{2})_{\tX^{-1}_1}\cdot \tX^2_1 =
(\tX^{-2}_{1})_{\tX_2} \cdot \tX^2_1 = c(\tX^{-
2}_{1})_{\tX^{-1}_2} \cdot \tX^2_1.$

On the other hand, $u^{-1}u_{\tX^{-1}_2} = \tX^2_2
\cdot (\tX^{-2}_{1})_{\tX^{-1}_2} \cdot \tX^{-2}_2 \cdot \tX^2_1 =
c(\tX^{-
2}_{1})_{\tX^{-1}_2} \cdot \tX^2_1.$ We get
$$u_{\tX^{-1}_2\tX^{-1}_1} = u^{-1}u_{\tX^{-1}_2}.\qquad\quad$$

$(2_b)$  Using (1), we can write
$$ u_{\tX_1\tX^{-1}_2} = (u^{-1}c)_{\tX^{-1}_2} = u^{-1}_{\tX^{-1}_2}
\cdot c.$$
Conjugating by $\tX_1\1$ and using $(2_a)$ we get, 
$$u_{\tX_1\tX^{-1}_2\tX^{-1}_1} =
(u^{-1})_{\tX^{-1}_2\tX^{-1}_1}\cdot c = u^{- 1}_{\tX^{-1}_2} u \cdot
c.$$
 Using (1) again we write 
$$u^{-1}_{\tX_1}
u_{\tX_1\tX^{-1}_2} = cu \cdot u^{-1}_{\tX^{-1}_2} c =
u^{-1}_{\tX^{-1}_2} u \cdot
c.$$ (We use $u_{X\1_2}=X_2^{-2}\cdot X_1^2$ and $[u,u_{\tX^{-1}_2}] =
[(\tX^2_{1})_{\tX^{-1}_2} \cdot \tX^{-2}_2, \tX^{-2}_2\tX^2_1] = c \cdot
c \cdot c = c)$.  
Thus, $$u_{\tX_1\tX^{-1}_2\tX^{-1}_1} = u^{-1}_{\tX_1}
u_{\tX_1\tX^{-1}_2}.\qquad$$

(3)  Clearly, $\forall j \geq 4$, $u_{\tX_j} = u$. 
We need to check that $u_{\tX_3}=u.$ Consider
$u_{\tX^{-1}_3} = (\tX^2_{1})_{\tX^{-1}_2\tX^{-1}_3}
(\tX^{-2}_{2})_{\tX^{-1}_3}$.

Since $u$ can also be written as $u = (\tX^2_{1})_{\tX^{-1}_2} \cdot
\tX^{-2}_2
= \tX^{-2}_2 \cdot (\tX^2_{1})_{\tX^{-1}_2} \cdot c = \tX^{-2}_2 \cdot
(\tX^2_{1})_{\tX_2}$, we have:  $$\align u_{\tX^{-1}_3} &= u
\Leftrightarrow
(\tX^2_{1})_{\tX^{-1}_2\tX^{-1}_3} \cdot( \tX^{-2}_{2})_{\tX^{-1}_3} =
\tX^{-2}_2 \cdot( \tX^2_{1})_{\tX_2} \Leftrightarrow \tX^2_2 \cdot
(\tX^2_{1})_{\tX^{-1}_2\tX^{-1}_3}\\& = (\tX^2_{1})_{\tX_2} \cdot
(\tX^2_{2})_{\tX^{-1}_3}.\endalign$$
The last equality is valid since $\{(\tX_{1})_{\tX_2},
\tX_2, (\tX_{2})_{\tX^{-1}_3}, (\tX_{1})_{\tX^{-1}_2\tX^{-1}_3}\}$ form
a good
quadrangle.
See Fig.   6 for the curves representing the quadrangle.
(Recall that $(X_1)_{X_2}=H((x_1)X_2)).$   \quad\qed\edm

\midspace{2.00in}
\caption{
Fig.   6}

\subheading{Construction of $\bold{\underline G(n)}$}

For $n \geq 3$ we define the group $\underline{G}(n)$ as follows:

Generators:  $s_1,u_1,u_2, \ldots, u_{n-1}$.

Relations:  $$\align&[s_1,u_i] = 1\quad \forall i = 1,3, \ldots
n-1;\quad
[u_i,u_j] = 1\ \text{when}\ |i-j| \geq 2;\\
&[s_1,u_2] = [u_i,u_{i+1}]   \quad \forall i = 1,2,
\ldots n - 2;\\
&[u_1,u_2] = [u_1,u_2]_{s_1} = [u_1,u_2]_{u_i} \;\quad \forall i =
1,2, \ldots, n - 1;\\
&[u_1,u_2]^2 = 1.\endalign$$
 
We denote: $[u_1,u_2]$ by $\nu.$
Clearly, $\nu\in\Center(\underline G(n)),$\ $\underline G(n)=\{1,\nu\}$
and $\nu^2=1.$

\bigskip

\subheading{Equivalent construction of $\bold{\underline G(n)}$}

Consider a free abelian group $A(n)$ with generators 
$S_1,V_1, \ldots, V_{n-1}$ and a skew-symmetric ${\Bbb Z}/2$-
valued bilinear form $Q(x,y)$ on $A(n)$ defined by: $Q(S_1,V_i) =
0 \; \forall i = 1,3, \ldots, n - 1;$ $Q(V_i,V_j) = 0$ when $|i-
j| \geq 2$, $Q(S_1,V_2) = Q(V_i,V_{i+1}) = 1$ $\forall i =
1,2,\ldots, n -2$. 
One can check that there exists a unique central extension $G$ of
$A(n)$ by ${\Bbb Z}/2$ with $\operatorname{Ab}(G) \simeq A(n)$\
($G^\prime
\simeq {\Bbb Z}/2)$ and  $[x,y] =
Q(\overline{x},\overline{y})$,  $\forall x,y \in G $ and
$\overline{x},$  $\overline{y}$ 
the images of $x,y$ in $A(n)$.
\proclaim{Claim   6.2}
\roster\item
The above  central extension is isomorphic to  $\underline{G}(n)$. 
\item $\Ab(\underline{G}(n))$ is a free abelian group with $n$
generators (i.e., $A(n)$) and $\underline G(n)^\prime \simeq {\Bbb Z}/2$
generated by $\nu=[u_1,u_2]$.  
\item
The following formulas define a
$\tilde{B}_n$-action on $\underline{G}(n)$ for  $\{\tX_1, \ldots,
\tX_{n-1}\},$   a frame of $\tilde{B}_n$.  $$\alignat3
&\underline{\tX_1\text{-action}}\quad &&
\underline{\tX_2\text{-action}}\quad && \underline{\tX_k\text{-action},
k
\geq 3} \\ 
&s_1 \rightarrow s_1 && s_1 \rightarrow u_2s_1&&
s_1 \rightarrow s_1 \\
&u_1 \rightarrow u^{-1}_1 \nu && u_1 \rightarrow
u_2u_1&& u_{k-1} \rightarrow u_ku_{k-1} 
\\& u_2
\rightarrow u_1u_2 && u_2 \rightarrow u^{-1}_2 \nu
&& u_k \rightarrow u^{-1}_k \nu\\
&u_j \rightarrow u_j \; \forall j \geq 3\quad && u_3 \rightarrow
u_2u_3; && u_{k+1} \rightarrow u_ku_{k+1}\\      
& &&u_j \rightarrow u_j \ \forall j \geq 4\quad && u_j
\rightarrow u_j \ \forall j \neq k-1,k,k+1
\endalignat$$
 \item
 Let $b \in \tilde{B}_n$, $y = (\tX_{1})_{b}$.  Then the $y^2$-action on
$\underline{G}(n)$ coincides with the conjugation by $(s_{1})_{b}$.
\item Let
$$ s_{ij}=\cases
s_1\
\qquad\quad\quad\qquad\qquad\qquad\quad\qquad\qquad\qquad\qquad\qquad\qquad
\text{if}\ (i,j)=(1,2); \\  (s_{1})_{\tX_2\dots \tX_{j-1}}\quad
\left(\overset\text{Claim   6.2(3)}\to=u_{j-1}\dots
u_2s_1\right)\qquad\quad\
\ \ \ \quad \quad\text{if}\ i=1,\ j\geq 2;\\ (s_1)_{\tX_2\dots
\tX_{j-1}\tX_1\dots
\tX_{i-1}}=\cases \nu\cdot u_{j-1}\dots u_1\cdot u_{i-1}\dots u_2s_1&\
\text{if}\ i\geq 3,\ j>i\\ \nu\cdot u_{j-1}\dots u_1s_1&\ \text{if}\
i=2,\
j>i.\endcases\endcases\tag1.20$$ 
Then   $$[s_{ij},s_{kl}] = 
\cases
\nu \quad& \text{if} \ \{i,j\} \cap \{k,l\} = 1;\\ 
1\quad& \text{otherwise}.
\endcases$$
\item Let $\tilde F_{n-1}$ be the subgroup
  of $\underline{G}(n)$ generated by $\{s_{n-
1,n},s_{n-2,n}, \ldots, s_{1n}\}$. 
qThe generators $\tX_1,\dots,\tX_{n-1}$ of  $\tB_{n-1}$ act on $\tilde
F_{n-1}$ as
follows: $$\align   (s_{jn})_{\tX_k} &=
s_{jn} \qquad   j \neq k, k + 1; \quad k=1 \ldots n - 2\tag1.21\\
(s_{k,n})_{\tX_k} &= s_{k+1,n} \\
(s_{k+1,n})_{\tX_k}&= s_{kn} \nu = s_{k+1,n}
s_{kn}s^{-1}_{k+1,n}.\endalign$$
 Their actions correspond to the
standard {\it Hurwitz moves} on ($s_{n-1,n},s_{n-1,n}, \ldots,
s_{1n})$ (see for example \cite{MoTe1}, Chapter $4$),  and it defines a
$\tB_{n-1}$-action.
\item There is a natural chain of embeddings $\underline{G}(3)
\subset \underline{G}(4) \subset \cdots \subset \underline{G}(n-
1) \subset \underline{G}(n)$ corresponding to the chain
$(s_1,u_1,u_2) \subset (s_1,u_1,u_2,u_3) \subset \ldots \subset
(s_1,u_1, \ldots u_{n-1})$.\endroster
\endproclaim

\pf 

(1), (2), and (3) are easy to verify.

(4)
Consider first the case $b =$ Id.  From (3) we get for the
$\tX^2_1$-action:
$$s_1 \rightarrow s_1, \; u_1 \rightarrow u_1, \; u_2 \rightarrow
u_2\nu, \; u_j \rightarrow u_j \; \forall j \geq 3.$$
At the same time by the first construction
$$(s_{1})_{s_1} = s_1,\ (u_1)_{s_1} = u_1,\; (u_2)_{s_1} = s^{-
1}_1u_2s_1 = u_2 \nu,\; (u_{j})_{s_1} = u_j \; \forall j \geq 3.$$
Thus $\tX^2_1$-action and $s_1$-conjugation coincide.  Consider now
any $b \in \tilde{B}_n$ and any $g \in \underline{G}(n)$.  Let $h
= g_{b^{-1}}$.  We have
$$g_{(\tX^2_{1})_{b}} = g_{b^{-1}\tX^2_1b} = ((h)_{\tX^2_1})_b =
(h_{s_1})_b = (h_{b})_{(s_{1})_{b}}=g_{(s_{1})_{b}}.$$  

(5), (6), (7) are easy to verify. \quad $\square$

\edm

\proclaim{Lemma   6.3}
Let $n \geq 3$.

Let $\{X_1,\dots,X_{n-1}\}$ be a frame of $B_n.$

Let
$$Z_{ij}=\cases X_1 & \text{if}\ (i,j)=(1,2);\\
(X_{1})_{X_2\dots X_{j-1}}& \text{if}\ i=1,\ j\geq 3;\\
(X_{1})_{X_2\dots X_{j-1}X_1\dots X_{i-1}}\quad&\text{if}\ i\geq
2,\ j>i.\endcases$$

Let $\tZ_{ij}$ be the image of $Z_{ij}$ in $\tilde B_n.$

Consider $\underline G(n)$ as a $ \tB_n$-group as in Claim \rom{  6.2}.

Then there exists a unique $\tilde{B}_n$-surjection $\Lambda_n:
\tilde{P}_n
\rightarrow \underline{G}(n)$ with $\Lambda_n(\tX^2_1) = s_1$
and $\Lambda_n(\tZ^2_{ij}) = s_{ij}$ for $1 \leq i < j \leq
n$. 
(In particular, $\Lambda_n(\tX_1^2)=s_1).$\ep

\demo{Proof}
Use induction on $n$.

 By definition $\tilde{P}_3$
 is obtained   from $P_3$ by adding the relations
$[\tZ^2_{12},\tZ^2_{23}] = [\tZ^2_{12},\tZ_{13}^2] = [\tZ^2_{23},
\tZ^2_{13}]$
and $[\tZ^2_{12},\tZ^2_{23}]^2 = 1$.
By Proposition 5.2 the embedding $P_3\subset P_4$ induces a homomorphism
$\tilde P_3\overset
i_2\to \rightarrow \tilde P_4.$ For $n = 3,  $ $\Lambda_3\: \tilde{P}_3
\rightarrow
\underline{G}(3)$ must be defined by $\lambda_3(\tZ^2_{ij}) =
s_{ij}$, $1 \leq i < j \leq 3$.  One can check directly that
$\Lambda_3$ is well
defined, and that it is a $\tilde{B}_3$-surjection.  Uniqueness of
such $\Lambda_3$ is evident.

Assume now that $n \geq 4$ and   the desired $\Lambda_{n-1}:
\tilde{P}_{n-1} \rightarrow \underline{G}(n-1)$ exists.

We shall establish first a surjection $\hat\Lambda_n:
P_n\twoheadrightarrow \underline G(n).$
Considering $\{X_1,X_2\} \subset \{X_1, X_2, X_3\} \subset \ldots
\subset \{X_1, \ldots, X_{n-1}\}$, we get a chain of embeddings
$B_3 \subset B_4 \subset \ldots \subset B_n$ and the
corresponding chain $P_3 \subset P_4 \subset \ldots \subset P_n$. 
To the latter corresponds a chain of homomorphisms:  $\tilde{P}_3
\overset{i_3}\to{\longrightarrow} \tilde{P}_4
\overset{i_4}\to{\longrightarrow} \ldots \rightarrow \tilde{P}_{n-1}
\overset{i_{n-1}}\to{\longrightarrow} \tilde{P}_n.$ 

The $Z_{ij}$'s defined here are in fact the same as those defined in
Remark 6.0 (by Claim
1.0(iii) and thus they generate $P_n$. 
For our proof it is better to use  these definitions of $Z_{ij}$'s.
 It is also known
that  $P_n \simeq P_{n-1} \ltimes F_{n-1},$ where $P_{n-1}$ is the
subgroup of $P_n$ generated
by $\{Z^2_{ij}, 1 \leq i < j \leq n - 1\}$, $F_{n-1}$ is the free
subgroup of $P_n$ generated by $\{Z^2_{in}, 1 \leq i \leq n - 1
\}$, and the semi-direct product $P_{n-1} \ltimes
F_{n-1}$ is defined according to the $P_{n-1}$-action on
$F_{n-1}$ which comes from the $B_{n-1}$-action by conjugation
(using $B_{n-1} \subset B_n \supset P_n$).  The latter coincides
with the standard $B_{n-1}$-action on $F_{n-1}$ (the generators
$X_{n-2}, \ldots, X_1$ of $B_{n-1}$
correspond to standard {\it Hurwitz moves} on $(Z^2_{n-1,n},
Z^2_{n-2,n}, \ldots , Z^2_{1n})$.  

  Using the
canonical map $P_{n-1} \rightarrow \tilde{P}_{n-1}$, we obtain from
$\Lambda_{n-1}$\ a $B_{n-1}$-surjection $\hat{\Lambda}_{n-1}:
P_{n-1} \rightarrow \underline{G}(n-1)$.  For the free subgroup
$F_{n-1}$ of $P_n$ generated by $\{Z^2_{in}, \; n - 1 \leq i \leq
1\}$ define $\mu_{n-1}: F_{n-1} \rightarrow \underline{G}(n)$ by
$\mu_{n-1}(Z^2_{in}) = s_{in}$.  Considering $P_n$ as $P_{n-1}
\ltimes F_{n-1}$, we  define
$\hat{\Lambda}_n: P_n
\rightarrow \underline{G}(n)$, which on $P_{n-1}$  
coincides with
$\hat{\Lambda}_{n-1}: P_{n-1} \rightarrow \underline{G}(n-1)
\subset \underline{G}(n)$ (see Claim   6.2(7)) and on $F_{n-1}$
coincides with
$\mu_{n-1}: F_{n-1} \rightarrow \underline{G}(n)$.  To show that such
$\hat{\Lambda}_n$ exists one has to check the following:

1)  The conjugation of $\mu_{n-1}(F_{n-1})$ by elements of
$\hat{\Lambda}_{n-1}(P_{n-1})(\subset \underline{G}(n))$
coincides with the $P_{n-1}$-action defined by $P_{n-1} \subset 
B_{n-1} \subset   B_n \rightarrow \tilde{B}_n$ and the given
$\tilde{B}_n$-action on $\underline{G}(n)$.  That
is, $\forall f \in \mu_{n-1}(F_{n-1})$ and $\forall h$ of the form $
\hat{\Lambda}_{n-1}(\tY)$ ($Y \in P_{n-1}$) we must have $h^{-1}fh
= f_{\tY}$.

2)  The $P_{n-1}$-action on $\mu_{n-1}(F_{n-1})$ (defined by
$P_{n-1} \subset B_{n-1} \subset B_n \rightarrow \tilde{B}_n$ and
the given $\tilde{B}_n$-action on $\underline{G}(n)$) comes from the
$B_{n-1}$-action on $\mu_{n-1}(F_{n-1})$ in which $X_{n-2},
\ldots X_1$ correspond to the standard {\it Hurwitz moves} on
$(s_{n-1,n},s_{n-2,n}, \ldots s_{1n})$.\edm

\demo{Proof of 1)}
Since $\forall b \in B_{n-1},$\   $\hat{\Lambda}_{n-1}((X^2_{1})_{b})
= (s_{1})_{b}$ we see from Claim   6.2 that $\forall f \in
\mu_{n-1}(F_{n-
1}),$\quad   $f_{\hat{\Lambda}_{n-1}((X^2_{1})_{b})} = (f_{s_1})_{b} =
f_{(\tX^2_{1})_{b}}$.  Since $P_{n-1}$ is generated by $\{Z^2_{ij}, 1
\leq i < j \leq n - 1\}$, i.e., by $\{(X^2_{1})_{b}, b \in B_n\}$
we get (1).\edm

\demo{Proof of 2)} It follows immediately from Claim   6.2(6).  

Thus, 1) and 2) are true and we can extend $\hat{\Lambda}_{n-1}$,
$\mu_{n-1}$ to a homomorphism $\hat{\Lambda}_n: P_n(= P_{n-1}
\ltimes F_{n-1}) \rightarrow \underline{G}(n)$. For  $1 \leq i < j \leq
n - 1$\quad
$\hat{\Lambda}_n(Z^2_{ij}) = \hat{\Lambda}_{n-1}(Z^2_{ij}) = s_{ij}$,
and for $1 \leq i \leq n
- 1$\quad $\hat{\Lambda}_n(Z^2_{in}) = \mu_{n-1}(Z^2_{in}) = s_{in}$.
Thus $\hat{\Lambda}_n(Z^2_{ij}) = s_{ij}$ for $1 \leq i < j
\leq n$.

Using induction, one can check directly that $\hat{\Lambda}_n$ is
a $B_n$-homomorphism (recall that by  Claim   6.2 we have explicit
formulas
for $s_{ij}$'s).

Since $s_{1n} = u_{n-1} \cdot\dots\cdot u_2s_1$ (see 1.20),  
$u_{n-1} = s_{1n}(u_{n-2} \ldots u_2s_1)^{-1}.$ Thus
$\underline{G}(n)$ is generated by $\underline{G}(n-1)$ and $s_{1n}.$
But $\underline G(n-1) =
\hat{\Lambda}_{n-1}(P_{n-1})$ $= \hat{\Lambda}_{n}(P_{n-1})$ and
$s_{1n} = \hat{\Lambda}_n(Z^2_{1n})$, and thus
$\hat{\Lambda}_n$ is a $B_n$-surjection.

Let ${N} = \ker (P_n \rightarrow \tilde{P}_n)\ (=\ker(B_n
\rightarrow \tilde{B}_n)$). Let $T = X^2_1X^2_3X^{- 2}_2Z^{-2}_{14}$. 
By Claim 1.2(b),
$ {N}$ is generated by $\{T_b,b \in B_n\}$.  We have
$\hat{\Lambda}_n(T) = \hat{\Lambda}_4(T) = s_1 \cdot s_{34}
\cdot s^{-1}_{23} s^{-1}_{14} \overset\text{Claim   6.2}\to= s_1 \cdot
\nu
u_3u_2 u_1 \cdot u_2s_1 \cdot s^{-1}_1u^{-1}_1 u^{-1}_2 \nu
\cdot s^{-1}_1u^{-1}_2 u^{-1}_3 = s_1 \nu u_3u_2\nu \cdot s^{-1}_1
u^{-1}_2 u^{-1}_3 =\nu u_3u_2u_2\1u_3\1$ (since $s_1$ commutes with
$u_3,$\
$[s_1,u_2]=\nu\in\Center\underline G(n)$ and $\nu^2=1)=\nu^2$ (since
$[u_3,u_2]=\nu)=$ Id. 
Since $\hat{\Lambda}_n$ is a $B_n$-homomorphism, we get
$\hat{\Lambda}_n(T_b) =$ Id $\forall b
\in B_n$, and thus $\hat{\Lambda}_n( {N}) = $ Id.  Hence
$\hat{\Lambda}_n$ defines canonically
a $\tilde{B}_n$-surjection $\Lambda_n: \tilde{P}_n \rightarrow
\underline{G}(n)$ with
$\Lambda_n(\tX^2_1) = s_1$.  

Uniqueness of such $\Lambda_n$ follows from the fact that
$\tilde{P}_n$ is generated by the $\tilde B_n$-orbit of $\tX^2_1$.
\qed\edm
 \proclaim{Theorem   6.4}
There exists a unique $\tilde{B}_n$-isomorphism $\Lambda_n:
\tilde{P}_n \rightarrow \underline{G}(n)$  with
$\Lambda_n(c) = \nu$ (see 5.3 and 6.2), s.t.  $\tilde{P}_{n,0}$ is
$\tilde{B}_n$-isomorphic to the subgroup
$G_0(n)$ of $\underline{G}(n)$, generated by $u_1, \ldots, u_{n-
1}$.  In
particular: 
\roster\item $\Ab\tilde{P}_n$ is a free abelian group with $n$
generators.
\item $\Ab\tilde{P}_{n,0}$ is a free abelian group 
generated by $n-1$ prime elements where $c$ is the commutator of any 2
of them.  
$\tilde P_n'=\tilde{P}^\prime_{n,0} \simeq {\Bbb Z}/2,$ generated by
$c.$
\item  $\tilde{P}_{n,0}$ is a primitive $\tilde{B}_n$-group
generated by the $\tilde{B}_n$-orbit of a prime element
${u} = \tX^2_1\tX_2^{-2}$, where $\tX_1,\tX_2$ are consecutive
half-twists
in $\tilde{B}_n$, $\tilde T = \tX_2\1\tX_1 \tX_2$ is the supporting
half-twist for
${u},$ and $c=[\tX_1^2,\tX_2^2]$ is its corresponding  central
element.\endroster
\ep

\demo{Proof}
 We first prove (3), since we use it to prove that the surjection
$\Lambda_n$ from Lemma 6.3
is in fact an isomorphism.
Complete $\tX_1,\tX_2$ to a frame $\tX_1,\dots,\tX_{n-1}$ of $\tilde
B_n.$
We apply conjugation by $\tX_2$ on  Lemma 6.1 and conclude using Lemma
2.2 that $u$ is prime,
its supporting half-twist is $\tX_2\1\tX_1\tX_2$ and its corresponding
central element is $c.$
By Remark 6.0,
$\tilde{P}_{n,0}$ is generated by $\{\tX^2_1\tZ^{-2}_{ij}, l \leq i < j
\leq n\}$.  Since
$\tX^2_1 \cdot \tZ^{-2}_{ij} = \tX^2_1 \tZ^{-2}_{1i} \cdot \tZ^2_{1i}
\cdot \tZ^{-2}_{ij}$ and
both $\tX^2_1\tZ^{- 2}_{1i}$,\ $\tZ^2_{1i} \cdot \tZ^{-2}_{ij}$ are
conjugates of $u$ (by
Claim 1.1(f)),   $\tilde{P}_{n,0}$ is generated by the $\tilde{B}_n$-
orbit of $u$.
Therefore,  $\tilde{P}_{n,0}$ is a
primitive $\tilde{B}_n$-group. 
Thus, we proved (3).

Polarize each $X_i$ (and $\tX_i$) according to the sequence $(X_1,
\ldots X_{n-1}$) (the ``end'' of $X_i$ = the ``origin'' of
$X_{i+1}$).  
By Theorem   3.3 $\forall i = 1, \ldots, n - 1$, $\exists$ a
unique prime element $\xi_i = L_{(u,\tX_1)}(\tX_i) \in
\tilde{P}_{n,0}$ such that $(\xi_i,\tX_i)$ is coherent
with
$(u,\tX_1)$.  Clearly $\xi_1 = u$.
By Lemma 3.1 the corresponding central element of $\xi_i,$\
$i=1,\dots,n,$ is also $c.$

By Proposition   4.1(1) and (2) we have $\forall i = 1, \ldots, n - 1$:
$$
(\xi_{i})_{\tX^{-1}_i} = \xi^{-1}_ic \qquad
(\xi_{i})_{\tX^{-1}_{i-1}} =
\xi_i \xi_{i-1}\qquad (\xi_{i})_{\tX^{-1}_{i+1}} = \xi_i\xi_{i+1}.
\tag1.22$$
It is also clear (Axiom (3)) that $\forall j \neq i$, $i - 1, i + 1$
$$
(\xi_{i})_{\tX_j} = \xi_i.
\tag1.23$$
We see from (1.22), (1.23) that the subgroup of $\tilde{P}_{n,0}$
generated by $\{\xi_1, \ldots, \xi_{n-1}\}$ is
closed under the
$\tilde{B}_n$-action.  Since $\tilde{P}_{n,0}$ is generated by
the $\tilde{B}_n$-orbit of $u = \xi_1$, we conclude
that
$\tilde{P}_{n,0}$ is generated by $\{\xi_1, \ldots,
\xi_{n-1}\}$.  This implies (see Remark 6.0) that $\tilde{P}_n$ is
generated by $\{\tX^2_1,\xi_1,\xi_2, \ldots,
\xi_{n-1}\}$.  

Since $X_1$ and $X_2$ are consecutive, $\tX_2=(\tX_1)_{\tX_2\1\tX_1\1}$
(Claim 1.0(ii)).
By the {\it uniqueness} of Theorem 3.3, $\xi_2\
(= L_{(\xi_1,\tX_1)}(\tX_2)) = (\xi_{1})_{\tX^{-1}_2\tX^{-1}_1}$ which
equals $
(\tX^{-2}_{1})_{\tX_2} \cdot \tX^2_1$\quad $(\xi_1=u).$ Thus 
$$
[\tX^2_1,\xi_2] = c.
\tag1.24$$

By Lemma   5.1 we have
$$
[\xi_i,\xi_j] =
\cases
c\quad & \text{if}\ |i-j| = 1 \\
1\quad &\text{if} \ |i-j| \geq 2.
\endcases\tag1.25
$$

Observe also that 
$$
(\tX^2_{1})_{\tX_2} = \xi_2 \cdot \tX^2_1.
\tag1.26$$

Formulas (1.22)--(1.26) show that we can define a
$\tilde{B}_n$-homomorphism
$M_n: \underline{G}(n) \rightarrow \tilde{P}_n$ with $M_n(s_1) =
\tX^2_1$,
$M_n(u_i) = \xi_i$, $i = 1,\ldots, n - 1$.  (See Claim   6.2.)

Since  $\tilde{P}_n$ is generated by the
$\tilde{B}_n$-orbit of $\tX^2_1$ and $\underline{G}(n)$ is generated by
the $\tilde{B}_n$-orbit of $s_1$, we conclude that $\Lambda_n$
and $M_n$ are inverses of each other.\quad $\square$\edm

The next corollary immediately follows from Theorem 6.4
\proclaim{Corollary 6.5}
\roster\item
We have the following sequence for $\tilde B_n:$

$1 <(\tilde P_n'=)\tilde P_{n,0}'<\tilde P_{n,0}<\tilde P_n<\tilde B_n$

where: 

$\tilde B_n/\tilde P_n\simeq S_n$

$\tilde P_n/\tilde P_{n,0}\simeq \Bbb Z$

$\tilde P_{n,0}/\tilde P_{n,0}'\simeq \Bbb Z^{n-1}$

$\tilde P_{n,0}'(=\tilde P_n')\simeq \Bbb Z_2$
\item
$\tilde B_n$ is an extension of a solvable group by a symmetric group.
In particualr, it is ``almost solvable'', i.e., it contains a solvable
group with a finite
index.\endroster\ep \bk

\subheading{\S7. Criterion for prime element}

The criterion for an element of a $\tilde B_n$-group to be prime,
presented in this section
(Proposition 7.1), is simpler to use than previous ones but its proof is
much longer.
The proof contains 8 lemmas.
This criterion will be used in the application of this paper  to study
fundamental groups
which turn out to be $\tilde B_n$-groups.

\proclaim{Proposition   7.1}
Assume $n \geq 5$.  Let $G$ be a $\tilde{B}_n$-group,
$(\tX_1,\tX_2, \ldots, \tX_{n-1})$ be a frame of $\tilde{B}_n$. 
Let $S$ be an element of $G$ with the following properties:
\roster\item"(0)"  $G$ is generated by $\{S_b,b \in \tilde{B}_n\}$;
\item"$(1_a)$"    $S_{\tX^{-1}_2\tX^{-1}_1} = S^{-1}S_{\tX^{-1}_2}$;
\item"$(1_b)$" $S_{\tX_1\tX^{-1}_2\tX^{-1}_1}$ = $S^{-
1}_{\tX_1}S_{\tX_1\tX^{-1}_2};$
\item"(2)"  For $\tau = SS_{\tX^{-1}_1},$ $T = S_{\tX^{-1}_2}$ we have:
\item"$(2_a)$"  $\tau_{\tX^2_1} = \tau$; 
\item"$(2_b)$" $\tau_T = \tau^{-1}_{\tX_1}$;
\item"(3)" $S_{\tX_j} = S$  $\forall j \geq 3$;
\item"(4)"  $S_c = S$, where $c = [\tX^2_1,\tX^2_2]$.\endroster

Then $S$ is a prime element of $G$, $\tX_1$ is the  supporting
half-twist of
$S$ and $\tau$ is the corresponding central element.  In particular,
$\tau^2 =
1$, $\tau \in\Center(G)$, $\tau_b = \tau$ $\forall$ $b \in \tilde{B}_n$.
\ep

\pf
The proof includes  several lemmas. From Theorem   5.2,
$c\in\Center(\tB_n),$\  $c^2 = 1.$  From Theorem   6.4 it follows that
$\tilde{P}^\prime_n$ is generated by $c$.  $\tilde{P}^\prime_n$ is a
normal
subgroup of $\tilde{B}_n$.  Write $\tilde{\tilde{B}}_n =
\tilde{B}_n/\tilde{P}^\prime_n$, $\tilde{\tilde{P}}_n =
\tilde{P}_n/\tilde{P}^\prime_n = \Ab \tilde{P}_n$. Clearly,
$\Tilde{\Tilde P}_n$ is a
commutative group.  We have $\Tilde{\Tilde\psi}_n:\Tilde{\Tilde
B}_n\twoheadrightarrow S_n.$
By abuse of notation we use $\psi_n$ for $\Tilde{\Tilde\psi}_n.$ Let
$Y\in B_n.$ 
By abuse of notation we denote the image of $Y$ in $\tB_n$ or in
$\Tilde{\Tilde
B}_n$ or in $\Tilde{\Tilde P}_n$  by the same symbol $\tY.$  It is clear
that
$\tilde{B}_n$ acts on $\Tilde{\Tilde{P}}_n$ (through conjugations) as
the
symmetric group $S_n = \tilde{B}_n/\tilde{P}_n$.

Since $S_c = S$ and $c\in\Center\tilde B_n,$ we have $\forall b \in
\tilde{B}_n$
$(S_b)_c = S_{bc} = (S_{c})_{b} = S_b$. Since $G$ is generated by
$\{S_b,b
\in \tilde{B}_n\}$ we have $\forall g
\in G,$ \ $g_c = g.$   In particular, we conclude that $\tilde{B}_n$
acts on $G$
as its quotient $\Tilde{\Tilde{B}}_n$; in other words, $G$ is a
$\Tilde{\Tilde{B}}_n$-group.

  Let $(D,K)$ be a model for $B_n,$\
${K} = \{a_1, \ldots, a_n\}$,\ $B_n=B_n[D,K].$ Take any $a_{i_1},a_{i_2}
\in
{K}$.  Let $\gamma_1,\gamma_2$ be two different simple paths in
$D - ({K} - a_{i_1}-a_{i_2})$ connecting $a_{i_1}$ with $a_{i_2}$, let
${H}(\gamma_1),{H}(\gamma_2)$ be the half-twists
corresponding to
$\gamma_1,\gamma_2$, and let ${\tilde H}(\gamma_1),
\tilde{H}(\gamma_2)$ be the images of ${
H}(\gamma_1), {H}(\gamma_2)$ in $\tilde{\tilde{B}}_n$. \enddemo

\proclaim{Lemma 1} Let $\g_1,\g_2$ be $2$ simple paths in
$D-\{K-a_{i_{1}}-a_{i_{2}}\}$ connecting $a_{i_{1}}$ with $a_{i_{2}}.$
Then:
$\tilde{H}(\gamma_1)^2 = \tilde H(\gamma_2)^2$. \ep
\demo{Proof of Lemma 1}
Choose a frame of ${B}_n$ $(Y_1, \ldots,
Y_{n-1})$ such that $Y_1 =  {H}(\gamma_1).$ Let $b \in
B_n$ be such that $\gamma_2 = (\gamma_1)b$, that is, ${H}(\gamma_1)_b
= {H}(\gamma_2)$. 
Let $\tY_i$ be the image of $Y_i$ in $\Tilde{\Tilde B}_n.$

Let $\sigma_1$ be the image of $b$ in $S_n$.  Since $(a_i)b =
a_i$, $(a_j)b = a_j$, $\sigma_1 \in
\operatorname{Stab}(i)\cap \operatorname{Stab}(j)$ in $S_n$. The
subgroup of
$\Tilde{\Tilde{B}}_n$ generated by $\tY_3, \ldots, \tY_{n-1}$ is mapped
by
$\Tilde{\Tilde \psi}_n:\Tilde{\Tilde B}_n\ri S_n$ onto
$\operatorname{Stab}(i)\cap \operatorname{Stab}(j).$  Choose $\tilde
b_1$ in
this subgroup with its image in $S_n$ equal to $\sigma_1$.  Clearly,
$(\tY_{1})_{\tilde b_1} = \tY_1$.   Since the image of
$\tilde b^{-1}_1\tilde{b}$ in $S_n$ is equal to $\sigma^{-1}_1\sigma_1 =
$ Id,
we have $\tilde b^{-1}_1 \tilde{b} \in \Tilde{\Tilde{P}}_n.$ Since
$\Tilde{\Tilde P}_n$ is commutative when considering $\tY_1^2$ as an
element
of $\Tilde{\Tilde P}_n,$\  $(\tY^2_{1})_{b^{-1}_1\tilde{b}} = \tY^2_1$. 
Thus,
we have 
$$\tilde{H}(\gamma_2)^2 = \tilde H(\gamma_1)^2_{\tilde{b}} =
(\tY^2_{1})_{\tilde b_1\tilde b^{-1}_1\tilde{b}} =
(\tY^2_{1})_{\tilde b^{-1}_1\tilde{b}} = \tY^2_1 = \tilde H(\gamma_1)^2  
\quad\square \quad\text{ for Lemma 1}$$\edm

\demo{Definition}\ $\underline{f_{ij}}$

$\forall i,j \in (1, \ldots, n)$, $i \neq j$, we define $f_{ij}
\in \Tilde{\Tilde{P}}_n$ as follows:  Take any simple path
$\gamma$ in $D - ({K}-a_i-a_j)$ connecting $a_i$ with
$a_j$.  Let $f_{ij} = \tilde{H}(\gamma)^2$.  Lemma 1 shows that this
definition does not depend on the choice of $\gamma$. We choose for
$i<j:$
$$f_{ij}=\cases (\tX_i^2)_{X_{i+1}}\cdot\dots\cdot X_{j-1}\quad &
2<j-i\\
\tX_i^2 & i+1=j.\endcases$$  

It
is clear that for  $\sigma_1$  the image   in $S_n$ of  $ b \in
\tilde{B}_n$ we have
$$ \quad(f_{ij})_b = f_{(i)\sigma_1,(j)\sigma_1}.$$ 

It is clear from our choice of $\g$ for $f_{ij}= \tilde
H(\gamma)^2$  that $$ 
\psi_n(\tilde{H}(\gamma)) = (i,j).$$
\medskip

It will be convenient to use the following notation for  $ g
\in G$ and $b \in \Tilde{\Tilde{B}}_n:$
\demo{Notation} \ $\un{[g,b]}:$  

For $g\in G$, $ b\in\Tilde{\Tilde B}_n$ and the action of $\Tilde{\Tilde
B}_n$
on $G$, we denote $[g,b] = g \cdot g^{-1}_{b^{-1}}.$
\medskip
One can check that:$$\align&g_b=g\Leftrightarrow [g,b]=1\\
&[g,b]_z=g_z(g_z)_{b_z\1}\1\\
&[g^{-1},b] = [g,b]^{-1}_g\\  &[g,b^{-1}] =
[g,b]^{-1}_b\\&[g_1g_2,b] = [g_2,b]_{g^{-1}_1} \cdot [g_1,b]\\
&[g,b_1b_2] = [g,b_1] \cdot [g,b_2]_{b^{-1}_1}.\endalign$$

\demo{Notation}\ $\underline{Q_{b,l,m}}$

$\forall b \in \Tilde{\Tilde{B}}_n,\quad \forall l,m \in (1, \ldots,
n),$ $l
\neq m$, we denote $Q_{b,l,m} = [S_b,f^{-1}_{lm}].$
\edm

\proclaim{Lemma 2}

\roster\item"(i)"Let $b \in \Tilde{\Tilde{B}}_n$ be such that
$(\{1,2\})\psi_n(b) \cap \{l,m\} = \emptyset$.  Then $Q_{b,l,m} = 1$.
\item"(ii)" Let $Q = Q_{\operatorname{Id},1,3} = [S,f^{-1}_{13}]$.
Then
$Q_{\tX^{-1}_2} = Q$.\endroster
\ep

\demo{Proof of Lemma 2}

(i)\ 
Let $\{l_1,m_1\} = (\{l,m\})\psi_n(b)^{-1}$.  So 
$(f_{lm})_{b\1}=f_{l_{1}m_{1}}.$ We have $\{1,2\} \cap \{l_1,m_1\} =
\emptyset$,
that is, $3 \leq l_1   , 3\le m_1$.  By our choice, $f_{lm}$ is a
product of
$X_j$ for $j\ge 3.$
Thus, using property (3) of $S$  we get
$S_{f_{l_1,m_1}} = S.$  In other words,  $[S,f_{l_1,m_1}\1] = 1$.  We
get
$$(Q_{b,l,m})_{b^{-1}} = [S_b,f^{-1}_{lm}]_{b^{-1}} =
[S,(f_{lm}\1)_{b^{-1}}] = [S,f_{l_1,m_1}\1] = 1,$$
and so $Q_{b,l,m} = 1$.   \quad\qed \quad for Lemma 2(i)

(ii)\
From $S_{\tX^{-1}_2\tX^{-1}_1} = S^{-1}S_{\tX^{-1}_2}$ (assumption
$(1_a$)
 of the 
Proposition) it follows that $S_{\tX^{-1}_2} = S
S_{\tX^{-1}_2\tX^{-1}_1}$.  Applying $\tX^{-1}_3$, we get
$$
S_{\tX^{-1}_2\tX^{-1}_3} = SS_{\tX^{-1}_2\tX^{-1}_1\tX^{-1}_3}
\tag1.27$$
which, after applying $\tX^{-1}_2$, gives
$$S_{\tX^{-1}_2\tX^{-1}_3\tX^{-1}_2} = S_{\tX^{-1}_2}
S_{\tX^{-1}_2\tX^{-
1}_1\tX^{-1}_3\tX^{-1}_2}.$$
Since $S_{\tX^{-1}_2\tX^{-1}_3\tX^{-1}_2} =
S_{\tX^{-1}_3\tX^{-1}_2\tX^{-
1}_3} = S_{\tX^{-1}_2\tX^{-1}_3}$, we obtain $S_{\tX^{-1}_2\tX^{-1}_3}
=$\newline $
S_{\tX^{-1}_2}S_{\tX^{-1}_2\tX^{-1}_1\tX^{-1}_3\tX^{-1}_2}$, or
$$
S_{\tX^{-1}_2} = S_{\tX^{-1}_2\tX^{-1}_3}
S^{-1}_{\tX^{-1}_2\tX^{-1}_1\tX^{-
1}_3\tX^{-1}_2}
\tag 1.28$$

Let  $b_1=\tX^{-
1}_2\tX^{-1}_1\tX^{-1}_3\tX^{-1}_2$. Observing that
$(f_{13})_{\tX^{-1}_2} =
f_{12}$, we get\linebreak from (1.28): $Q_{\tX^{-1}_2} =
[S_{\tX^{-1}_2},(f^{-1}_{13})_{\tX^{-1}_2}] = [S_{\tX^{-
1}_2\tX^{-1}_3}S^{-1}_{b_1},f^{-1}_{12}].$  Thus
$$Q_{\tX^{-1}_2} 
 = [S^{-1}_{b_1},f^{-
1}_{12}]_{S^{-1}_{\tX^{-1}_2\tX^{-1}_3}} \cdot
[S_{\tX^{-1}_2\tX^{-1}_3},
f^{-1}_{12}].
\tag1.29$$

Since $\psi_n(b_1) = (2\quad 3)\ (1\quad 2)\ (3\quad 4)\ (2\quad 3)$
(products
of transpositions),\linebreak $(\{1,2\})\psi_n(b_n) = \{3,4\}.$ Since
$\{3,4\}
\cap \{1,2\} = \emptyset$, we get from   (i)
 that $Q_{b_1,1,2} = 1.$ Thus, $[S^{-1}_{b_1},f^{-1}_{12}] =
[S_{b_1},f^{-1}_{12}]^{- 1}_{S_{b_1}} = (Q^{-1}_{b_1,1,2})_{S_{b_1}}$.
(1.29) now gives
$$
Q_{\tX^{-1}_2} = [S_{\tX^{-1}_2\tX^{-1}_3},f^{-1}_{12}].
\tag1.30$$

Consider a quadrangle formed by $\{a_1,a_2,a_3,a_5\}$, as in
Fig.   7.1.  By Lemma   1.2, we can write in $\Tilde{\Tilde{P}}_n$:
$f_{35}f_{12} = f_{25}f_{13}$, or $f_{12} = f_{35}^{-
1}f_{25}f_{13}$, $f^{-1}_{12} = f^{-1}_{13} f^{-
1}_{25}f_{35}$.

\midspace{1.40in}
\caption{Fig.   7.1}

>From (1.30) we get
$$
Q_{\tX^{-1}_2} =
[S_{\tX^{-1}_2\tX^{-1}_3},f^{-1}_{13}f^{-1}_{25}f_{35}] =
[S_{\tX^{-1}_2\tX^{-1}_3},f^{-1}_{13}] \cdot [S_{\tX^{-1}_2\tX^{-
1}_3},f^{-1}_{25}f_{35}]_{f_{13}}
\tag1.31$$

Consider $[S_{\tX^{-1}_2\tX^{-1}_3},f^{-1}_{25}f_{35}] = [S_{\tX^{-
1}_2\tX^{-1}_3}, f^{-1}_{25}] \cdot [S_{\tX^{-1}_2\tX^{-
1}_3},f_{35}]_{f_{25}} = Q_{\tX^{-1}_2\tX^{-1}_3,2,5} \cdot [S_{\tX^{-
1}_2\tX^{-1}_3},f^{-1}_{35}]^{-1}_{f^{-1}_{35}f_{25}} = Q_{\tX^{-
1}_2\tX^{-1}_3,2,5} \cdot (Q^{-1}_{\tX^{-1}_2\tX^{-1}_3,3,5})_{f^{-
1}_{35}f_{25}}.$  Since $\psi(\tX^{-1}_2\tX^{-1}_3) =
(2\quad 3)(3\quad 4),$ the images of $\{1,2\}$ under it are  $
\{1,4\}.$ But $\{1,4\} \cap \{2,5\} = \emptyset$ and $\{1,4\} \cap
\{3,5\} =
\emptyset.$ Thus,  we get by (i) that
$Q_{\tX^{-1}_2\tX^{-1}_3,2,5} = Q_{\tX^{-1}_2\tX^{-1}_3,3,5} = 1$, and
so $[S_{\tX^{-1}_2\tX^{-1}_3},f^{-1}_{25}f_{35}] = 1$.  (1.31)
now implies $Q_{\tX^{-1}_2} = [S_{\tX^{-1}_2\tX^{-1}_3}, f^{-1}_{13}]$. 
By (1.27) $S_{\tX_2^{-1}\tX^{-1}_3} = S \cdot
S_{\tX^{-1}_2\tX^{-1}_1\tX^{-1}_3}$ which gives
$$\align Q_{\tX^{-1}_2} &= [S \cdot S_{\tX^{-1}_2\tX^{-1}_1\tX^{-1}_3},
f^{-
1}_{13}] = [S_{\tX^{-1}_2\tX^{-1}_1\tX^{-1}_3},f^{-1}_{13}]_{S^{-1}}
\cdot [S,f^{-1}_{13}] \\&=
(Q_{\tX^{-1}_2\tX^{-2}_1\tX^{-1}_3,1,3})_{S^{-1}} \cdot Q.\endalign$$

The value of $ \psi(\tX^{-1}_2\tX^{-1}_1\tX^{-1}_3)\ ( =
 (2 \quad 3)(1\quad 2)(3\quad 4))$ on $\{1,2\}$ is $\{2,4\}$.  Since
$\{2,4\}
\cap \{1,3\} = \emptyset$ we get from part (i) of the Lemma that
$Q_{\tX^{-1}_2\tX^{-1}_1\tX^{-1}_3,1,3} = 1$, therefore,
$$Q_{\tX^{-1}_2} = Q \qed\quad \text{for   Lemma 2(ii).}$$
\edm
\proclaim{Lemma 3} $\tau=Q\1.$\ep
\demo{Proof of Lemma 3}
By the assumption on $\tau,$\ $S_{\tX^{-1}_1} = S^{-1}
\tau.$  By definition of $T,$\  $T_{\tX^{-2}_1} =
S_{\tX^{-1}_2\tX^{-2}_1}.$
We apply assumption $(1_a)$ twice to get, using $S_{\tX_{1}\1}=S\1\tau,$
that
\newline  $(S^{- 1}S_{\tX^{-1}_2})_{\tX^{-1}_1} =
S^{-1}_{\tX^{-1}_1}S_{\tX^{-1}_2\tX^{- 1}_1} = \tau^{-1}S \cdot
S^{-1}S_{\tX^{-1}_2} = \tau^{-1} S_{\tX^{- 1}_2} = \tau^{-1}T$.  Thus
$$
T_{\tX^{-2}_1} = \tau^{-1}T,
$$
or
$$\tau^{-1} = T_{\tX^{-2}_1} T^{-1}.\tag1.32$$
Applying $\tX^2_1$ on (1.32) and using $\tau_{\tX^2_1} = \tau$
(assumption $(2_a)$), we get $$\align\tau^{-1}& = T \cdot
T^{-1}_{\tX^2_1} = [T,\tX^{- 2}_1] = [S_{\tX^{-1}_2},f^{-1}_{12}] =
[S,(f^{-1}_{12})_{{\tX_2}}]_{\tX^{- 1}_2} =
[S,f^{-1}_{13}]_{\tX^{-1}_2}\\& =
Q_{\tX^{-1}_2} \overset\text{by 
 Lemma 2}\to =
Q,\endalign$$ that is, $$
\tau^{-1}=Q, \quad \text{or} \quad \tau =
Q^{-1}.\quad\square\quad\text{for
Lemma 3} \tag 1.33$$

\proclaim{Lemma 4}
$\forall j \geq 3$, $\tau_{\tX_j} = \tau.$
\ep

\demo{Proof of Lemma 4}
From $\tau = SS_{\tX^{-1}_1}$ and $S_{\tX_j} = S$ $\forall j \geq 3$
it follows that $\tau_{\tX_j} = \tau$ $\forall j \geq 3$.    \qed\quad
for
Lemma 4.\edm

\proclaim{Lemma 5}
$\tau_{\tX_1} = \tau.$
\ep

\demo{Proof of Lemma 5}  Let us use
now $\tau^{-1}_{\tX_1} = \tau_T$ ($(2_b)$ of the Proposition).

By Lemmas 2 and 3  $\tau_{\tX_2} = \tau$. 

 Thus,
$\tau_T = \tau_{S_{\tX^{-1}_2}} = S^{-1}_{\tX^{-1}_2} \tau
S_{\tX^{-1}_2} = (S^{-1}\tau S)_{\tX^{-1}_2} = (\tau_S)_{\tX^{-
1}_2}$.

 So $\tau_{\tX^{-1}_1} = \tau_{\tX_1} =\tau_T\1= (\tau^{-
1}_S)_{\tX^{-1}_2}$, or 
$$
\tau_S = \tau^{-1}_{\tX^{-1}_1\tX_2}.
\tag1.34$$

Since $\tau_{\tX_3} = \tau$ and $S_{\tX_3} = S$, we get
$(\tau_S)_{\tX_3} = \tau_S$ and
$$
\tau_{\tX^{-1}_1\tX_2\tX_3} = (\tau^{-1}_S)_{\tX_3} = \tau^{-1}_S =
\tau_{\tX^{-1}_1\tX_2}.
\tag1.35$$

Applying $\tX^{-1}_2$ on (1.35), we get 
$
\tau_{\tX^{-1}_1\tX_2\tX_3\tX_2^{-1}} = \tau_{\tX^{-1}_1}.
 $
Since $\tau_{X_3}=\tau$ and $\la X_2,X_3\ra=1,$\ 
$\tau_{\tX^{-1}_1\tX_2\tX_3\tX_2^{-1}} = \tau_{\tX^{-1}_1\tX^{-
1}_3\tX_2\tX_3}
= \tau_{\tX^{-1}_3\tX^{-1}_1\tX_2\tX_3} = \tau_{\tX^{-1}_1\tX_2\tX_3}$.
Thus
$$\tau_{\tX_{1}\1\tX_{2}\tX_{3}}=\tau_{\tX_{1}\1}\tag1.36.$$
Combining
formulas (1.35)-(1.36) we get $\tau_{\tX^{-1}_1} =  
\tau_{\tX^{-1}_1\tX_2}$.  Applying it to $\tX_1$ we get $\tau=
\tau_{\tX^{-1}_1\tX_2\tX_1 } = \tau_{\tX_2\tX_1\tX^{-1}_2}  =
\tau_{\tX_1\tX^{-1}_2}$. Thus $\tau = \tau_{\tX_1\tX^{-1}_2}$,
or $\tau_{\tX_1} = \tau_{\tX_2} = \tau$. \linebreak 
\qed\quad for  Lemma 5.\edm

\proclaim{Lemma 6} 
$\tau_{\tX_j} = \tau$ $\forall j = 1,2,\ldots, n - 1$.
\ep
\demo{Proof of Lemma 6} By Lemmas 2, 3, 4, 5.\quad \qed\edm
\proclaim{Lemma 7} $\tau_S = \tau^{-1}$.
\ep
\demo{Proof of Lemma 7}
From $\tau_{\tX_1} = \tau_{\tX_2} = \tau$ and (1.34).\quad \qed
\edm

\proclaim{Lemma 8}
$\tau_S = \tau$.
\ep

\demo{Proof of Lemma 8}
Consider assumption $(1_b)$ of the Proposition
$$S_{\tX_1\tX^{-1}_2\tX^{-1}_1} = S^{-1}_{\tX_1}S_{\tX_1\tX^{-1}_2}.$$

Using $S_{\tX^{-1}_1} = S^{-1} \tau$ and $\tau_{\tX_1} = \tau$, we
get $S = S^{-1}_{\tX_1}\tau$, or $S^{-1}_{\tX_1} = S\tau^{-1}$,
$S_{\tX_1} = \tau S^{-1}$.  Assumption $(1_b)$ now gives (using
$\tau_{\tX_i}
= \tau \; \forall i = 1,\ldots, n - 1)$
$$\tau S^{-1}_{\tX^{-1}_2\tX^{-1}_1} = S\tau^{-1} \cdot \tau S^{-
1}_{\tX^{-1}_2} = SS^{-1}_{\tX^{-1}_2} = ST^{-1}.$$

On the other hand, by $(1_a)$ and (2), $S_{\tX^{-1}_2\tX^{-1}_1} =
S^{-1}T.$ 
Thus $$\tau S_{X_{2}\1}\1X_1\1=\tau T\1S.$$

We compare the last 2 expressions to get
$\tau T^{-1}S = ST^{-1}$, or
$$
\tau = ST^{-1}S^{-1}T, \quad\text{or} \quad T^{-1}_{S^{-1}} =
\tau T^{-1}.
\tag1.37$$

By Lemmas 3 and 6 $  Q  =
Q_{\tX^{-1}_1\tX^{-1}_2\tX^{-1}_3} =
[S_{\tX^{-1}_1\tX^{-1}_2\tX^{-1}_3},
(f^{-1}_{13})_{\tX^{-1}_1\tX^{-1}_2\tX^{-1}_3}].$
Thus
$$Q= [S_{\tX^{-1}_1\tX^{-1}_2\tX^{-1}_3}, f^{-1}_{24}]
\tag 1.38 $$
(we use $(\{1,3\})\psi(\tX^{-1}_1\tX^{-1}_2\tX^{-1}_3) =
(\{1,3\}) (1\quad 2)(2\quad 3)(3\quad 4) = \{4,2\})$.

Considering a quadrangle formed by $a_2,a_3,a_4,a_5$ (see Fig.   7.2)

\midspace{2.00in}
\caption{Fig.   7.2}

\flushpar we can write in $\Tilde{\Tilde{P}}_n$ (Lemma   2.2) 
$f_{35}f_{24} =
f_{25}f_{34}$, or $f_{24} = f^{-1}_{35}f_{25}f_{34}$, $f^{-
1}_{24} = f^{-1}_{34}f^{-1}_{25}f_{35}$.  From (1.38) we get,
denoting by $b = \tX^{-1}_1\tX^{-1}_2\tX^{-1}_3$,
$$\align Q &= [S_{\tX^{-1}_1\tX^{-1}_2\tX^{-1}_3},f^{-1}_{24}] = [S_b,
f^{-
1}_{34}f^{-1}_{25}f_{35}]\tag1.39\\
&= [S_b,f^{-1}_{34}][S_b,f^{-1}_{25}f_{35}]_{f_{34}} = Q_{b,3,4}
\cdot [S_b,f^{-1}_{25}]_{f_{34}} \cdot
[S_b,f_{35}]_{f_{25}f_{34}}\\
&
= Q_{b,3,4} \cdot (Q_{b,2,5})_{f_{34}} \cdot [S_b,f^{-1}_{35}]^{-
1}_{f^{-1}_{35}f_{25}f_{34}}= Q_{b,3,4} \cdot (Q_{b,2,5})_{f_{34}} \cdot
(Q_{b,3,5})_X\1. \endalign$$
Now,  $(\{1,2\})\psi(b) = \{1,2\}(1\quad 2)(2\quad 3)(3\quad 4) =
\{4,1\}$.  Since $\{4,1\} \cap \{2,5\} = \emptyset$ and
$\{4,1\} \cap \{3,5\} = \emptyset$, we get by   Lemma 2 that  
$Q_{b,2,5} = Q_{b,3,5} = 1$, and by (1.39)
$$Q = Q_{b,3,4}.$$

We can write $S_{\tX^{-1}_1\tX^{-1}_2\tX^{-1}_3} = (S^{-1}\tau)_{\tX^{-
1}_2\tX^{-1}_3} = S^{-1}_{\tX^{-1}_2\tX^{-1}_3} \tau = T^{-1}_{\tX^{-
1}_3}\tau$ (using $\tau_{\tX_i} = \tau$ $\forall i$).  So $Q =
Q_{b,3,4} = [S_{\tX^{-1}_1\tX^{-1}_2\tX^{-1}_3},f^{-1}_{34}] 
\overset\text{by}\ f_{34} = \tX^2_3\to = [T^{-1}_{\tX^{-1}_3} \tau,
\tX^{-2}_3] =
[\tau, \tX^{-2}_3]_{T_{\tX^{-1}_3}} \cdot [T^{-1}_{\tX^{-1}_3}, \tX^{-
2}_3] = [T^{-1}_{\tX^{-1}_3}, \tX^{-2}_3] =  [T^{-1},\tX^{-2}_3]_{\tX^{-
1}_3}$.  Since $Q_{\tX_3} = Q$, we get
$$
Q = [T^{-1},\tX^{-2}_3].
\tag1.40$$

This implies $Q_{S^{-1}} = (T^{-1}T_{\tX^2_3})_{S^{-1}}\overset\text{
assumption 3}\to = T^{- 1}_{S^{-1}} \cdot (T_{S^{-1}})_{\tX^2_3} =
[T^{-1}_{S^{-1}},\tX^{- 2}_3]$\linebreak $ \overset\text{by (1.37)}\to=
[\tau
T^{-1}, \tX^{-2}_3] = [T^{-1},\tX^{- 2}_3]_{\tau^{-1}} [\tau,\tX^{-2}_3]
\overset\text{by Lemma 4}\to= [T^{-1},\tX^{-2}_3]_{\tau^{- 1}}
\overset\text{by
(1.40)}\to= Q_{\tau^{-1}}.$

Using $Q=\tau\1$ we get $\tau_{S\1}\1=\tau_{\tau\1}\1.$
Thus, $\tau_{S^{-1}} = \tau$ and $\tau_S = \tau$.   \newline\qed\quad
for  Lemma
8.

We can now finish the proof of Proposition   7.1.  

By Lemma   2.4, we only have to prove that $\tau^2 = 1$,
$\tau_b = \tau$ $\forall b \in \Tilde{\Tilde{B}}_n$ and $\tau
\in\operatorname{Center}(G)$.

By the previous Lemma, $\tau_S = \tau$, and by Lemma 7, $\tau_S =
\tau^{-1}$. 
Thus, $\tau = \tau^{-1}$ and    $\tau^2 = 1$.

By Lemma 6, $\tau_{\tX_i} = \tau$
$\forall i \in (1, \ldots, n - 1)$.  Thus $\tau_b = \tau$ $\forall b \in
\Tilde{\Tilde{B}}_n$.

 By Lemma 8,
$\tau_S = \tau$, i.e. $[\tau,S]=1.$
Let $b\in\tB_n:$ \ $[\tau,S_b] = [\tau_{b^{-1}},S]_b =
[\tau,S]_b  = 1$.   

Thus $\tau$ commutes with $S_b$\ $\forall b\in \Tilde{\Tilde
B}_n.$ Since $
\ker(\tB_n\ri\Tilde{\Tilde B}_n)$ acts trivially on $G$,\linebreak
$\Tilde{\Tilde B}_n$   acts on $G$ via $\tilde B_n,$ and thus $\tau$
commutes
with $S_b$\ $\forall b\in \tB_n.$

 By assumption (0) of the proposition, $G$ is
generated by $\{S_b\}_{b\in\tB_n}.$ Thus   $\tau \in\Center(G).$ 

\hfill\qed \quad for Proposition   7.1\edm

\bk

\end
\Refs\widestnumber\key{MoTe10}
\ref\key A
\by  Artin E.\paper Theory of braids\jour Ann. Math. \vol 48\pages 
101-126\yr 1947
\endref
\ref\key B \by Birman J.\book Braids, Links and Mapping Class Groups. 
\publ Princeton University Press\yr 1975\endref


\ref\key MoTe1  \by Moishezon B., Teicher M. \paper Braid group
technique in complex geometry, I \jour Contemp. Math. \vol 78 \yr 1988
\pages 425-555\endref


\ref\key MoTe2 \by  Moishezon B., Teicher M. \paper
Fundamental groups of complements of curves in $\Bbb C\Bbb P^2$ as
solvable groups \jour IMCP \vol 9\yr 1996\pages 329-325\endref
 

\ref\key Te\by Teicher M.\paper Barid groups, algebraic surfaces and
fundamental groups of complements of curves\jour Proceedings of the
Santa
Cruz Summer Institute 1995 \endref\endRefs
